\documentclass{article}
\usepackage[a4paper, total={6in, 9in}]{geometry} 
\usepackage{graphicx} 
\usepackage{amsthm}
\usepackage{amsfonts}
\usepackage{amssymb}
\usepackage{amsmath}
\usepackage{mathtools}
\usepackage[maxbibnames=99]{biblatex}
\usepackage{xcolor}
\usepackage{esint} 
\usepackage{tikz-cd}
\usepackage{blindtext}
\usepackage{titlesec}
\usepackage[hyphenbreaks]{breakurl}

\usepackage[citecolor=blue,colorlinks]{hyperref}

\newtheorem{innercustomgeneric}{\customgenericname}
\providecommand{\customgenericname}{}
\newcommand{\newcustomtheorem}[2]{%
  \newenvironment{#1}[1]
  {%
   \renewcommand\customgenericname{#2}%
   \renewcommand\theinnercustomgeneric{##1}%
   \innercustomgeneric
  }
  {\endinnercustomgeneric}
}

\newcustomtheorem{customthm}{Theorem}
\newcustomtheorem{customC}{Corollary}
\newcustomtheorem{assumption}{Assumption}

\newtheorem{thm}{Theorem}[section]

\newtheorem{lemma}[thm]{Lemma}
\newtheorem{proposition}[thm]{Proposition}
\newtheorem*{corollary*}{Corollary}

\theoremstyle{definition}
\newtheorem{definition}[thm]{Definition}
\newtheorem{remark}[thm]{Remark}
\newtheorem{example}[thm]{Example}

\newcommand{\graph}{\mathsf{graph}}
\newcommand{\bb}[1]{\mathbb{#1}}
\newcommand{\ssf}[1]{\mathsf{#1}}
\newcommand{\supp}{\mathsf{supp}}

\newcommand{\lip}{\mathsf{lip}}
\newcommand{\Lip}{\mathsf{Lip}}

\newcommand{\Ch}{\mathsf{Ch}}
\newcommand{\Ric}{\mathsf{Ric}}

\newcommand{\RCD}{\mathsf{RCD}}

\newcommand{\aH}{\mathsf{H}}
\newcommand{\BV}{\mathsf{BV}}
\newcommand{\W}{\mathsf{W}}
\newcommand{\X}{\mathsf{X}}
\newcommand{\M}{\mathsf{M}}
\newcommand{\A}{\mathsf{A}}

\newcommand{\sP}{\mathsf{P}}
\newcommand{\Hyp}{\mathsf{Hyp}}
\newcommand{\sd}{\mathsf{d}}
\newcommand{\sL}{\mathsf{L}}
\newcommand{\m}{\mathfrak{m}}

\newcommand{\mms}[3]{(\mathsf{#1},\mathsf{#2},\mathfrak{#3})}
\newcommand{\pmms}[4]{(\mathsf{#1},\mathsf{#2},\mathfrak{#3},\mathsf{#4})}

\makeatletter
\newsavebox\myboxA
\newsavebox\myboxB
\newlength\mylenA

\newcommand*\xoverline[2][0.75]{%
    \sbox{\myboxA}{$\m@th#2$}%
    \setbox\myboxB\null
    \ht\myboxB=\ht\myboxA%
    \dp\myboxB=\dp\myboxA%
    \wd\myboxB=#1\wd\myboxA
    \sbox\myboxB{$\m@th\overline{\copy\myboxB}$}
    \setlength\mylenA{\the\wd\myboxA}
    \addtolength\mylenA{-\the\wd\myboxB}%
    \ifdim\wd\myboxB<\wd\myboxA%
       \rlap{\hskip 0.5\mylenA\usebox\myboxB}{\usebox\myboxA}%
    \else
        \hskip -0.5\mylenA\rlap{\usebox\myboxA}{\hskip 0.5\mylenA\usebox\myboxB}%
    \fi}
\makeatother

\title{Convergence of the area functional on spaces with lower Ricci bounds and applications}
\author{
\footnote{Mathematical Institute, University of Oxford, Oxford, UK. E-mail address: alessandro.cucinotta@maths.ox.ac.uk} Alessandro Cucinotta}

\begin{document}

\maketitle

\begin{abstract}
We show that the heat flow provides good approximation properties for the area functional on proper $\RCD(K,\infty)$ spaces, implying that in this setting the area formula for functions of bounded variation holds and that the area functional coincides with its relaxation. We then obtain partial regularity and uniqueness results for functions whose hypographs are perimeter minimizing.
Finally, we consider sequences of $\RCD(K,N)$ spaces and we show that, thanks to the previously obtained properties, Sobolev minimizers of the area functional in a limit space can be approximated with minimizers along the converging sequence of spaces.
Using this last result, we obtain applications on Ricci-limit spaces.
\end{abstract}

\tableofcontents

\section*{Introduction}

In recent years there has been a growing interest in the study of metric (measure) spaces with a notion of lower bound on the Ricci curvature.
Initially, this was motivated by the interest in Ricci limits, i.e. spaces arising as Gromov-Hausdorff limits of manifolds with a uniform lower Ricci curvature bound. 
Later on, intrinsic formulations of such curvature bound were found, leading to the theory of $\RCD$ spaces. \par 
An $\RCD(K,N)$ space is a metric measure space where $K \in \bb{R}$ plays the role of a lower bound on the Ricci curvature and $N \in [1,+\infty]$ plays the role of an upper bound on the dimension. 
This class includes measured Gromov-Hausdorff limits of smooth manifolds of fixed dimension with uniform Ricci curvature lower bounds and finite dimensional Alexandrov spaces with sectional curvature bounded from below. Moreover, the $\RCD(K,N)$ class is closed with respect to measured Gromov-Hausdorff convergence. For this reason, 
in the $\RCD$ setting, the study of objects that retain their properties along converging sequences of spaces is particularly relevant.
This is one of the motivations that led to an extensive study of sets of finite perimeter and perimeter minimizing sets in the non-smooth context of $\RCD(K,N)$ spaces. 
\par 
In the paper, we study some approximation properties of the area functional on $\RCD$ spaces and Ricci limits, extending a series of results that were previously obtained for sets of finite perimeter in this context. We then derive some new applications.
The first central result concerns the possibility of approximating functions via the heat flow on proper $\RCD(K,\infty)$ spaces in a way that involves the area. This is the content of Theorem \ref{CT0}.
Given a metric measure space $(\X,\sd,\m)$, an open set $\Omega \subset \X$ and $u \in \BV(\Omega)$, we denote by $|D u|$ its total variation measure, by $|\nabla u|$ the density of the absolutely continuous part of $|Du|$ w.r.t. $\m$ and by $|D^s u|$ the singular part of $|Du|$. For every measurable set $E \subset \Omega$, we then define the area of $u$ in $E$ as
\[
\A(u,E):=\int_E \sqrt{1+|\nabla u|^2} \, d\m +|D^s u|(E).
\]
We also define
$\Hyp(u):=\{(x,t) \subset \Omega \times \bb{R}:t \leq u(x)\}$ and we denote by $\sL^0(\Omega)$ the set of measurable functions on $\Omega$. Given two sets $A,B \subset \X$, we denote by $\sP(A,B)$ the perimeter of $A$ in $B$. Finally, we denote the heat flow on an $\RCD(K,\infty)$ space by $P_t(\cdot)$.

\begin{customthm}{1}
\label{CT0}
    Let $(\X,\sd,\m)$ be a proper $\RCD(K,\infty)$ space and let $u \in \sL^{\infty}(\X) \cap \BV(\X)$. Then, for every bounded open set $E \subset \X$ such that $|D^su|(\partial E)=0$,
it holds
\[
\lim_{t \to 0} \A(P_t(u),E) = \A(u,E).
\]
\end{customthm}

As a first application of Theorem \ref{CT0} we obtain that the area formula holds on proper $\RCD(K,\infty)$ spaces, i.e. the area of a function of bounded variation coincides with the perimeter of its hypograph. 
This is the content of Theorem \ref{CT1}. 
Before stating the theorem, we recall that related results were proved for $PI$ spaces in \cite{Tens}, where it was conjectured that the area formula for functions of bounded variation should hold in that setting. 
In particular, \cite[Theorem $5.1$]{Tens} implies that the area formula holds for Sobolev functions on $PI$ spaces where weak and relaxed gradients coincide (and finite dimensional $\RCD(K,N)$ spaces fit in this framework).
Such formula was then proved for $\BV$ functions on finite dimensional $\RCD(K,N)$ spaces in \cite{Brena22}, while in the next theorem we extend the result to the infinite dimensional framework of proper $\RCD(K,\infty)$ spaces.

\begin{customthm}{2} \label{CT1}
Let $(\X,\sd,\m)$ be a proper $\RCD(K,\infty)$ space, let $\Omega \subset \X$ be open and bounded and let $u \in \sL^1(\Omega)$. Then $\sP(\Hyp(u),\Omega \times \bb{R})<+\infty$ if and only if $u \in \BV(\Omega)$. In this case, for every Borel set $E \subset \Omega$, it holds
    \[
\sP(\Hyp(u),E \times \bb{R})=\A(u,E).
\]
Moreover, if $u \in \sL^0(B_{2r}(x))$ and $\sP(\Hyp(u),B_{2r}(x) \times \bb{R})<+\infty$, then $u \in \BV(B_r(x))$ and the previous equation holds for every Borel set $E \subset B_r(x)$.
\end{customthm}

As anticipated, the previous result complements the one proved for doubling spaces supporting a Poincaré inequality in \cite{Tens}, where it was shown that if $(\X,\sd,\m)$ is such a space and $u \in \BV(\Omega)$, then
\[
\sP(\Hyp(u),\Omega \times \bb{R}) \geq \A(u,\Omega).
\]
Note that, although $\RCD(K,\infty)$ spaces are not necessarily doubling, the proof of the inequality given in \cite{Tens} works with very slight changes also in our setting. On the other hand, the proof that we give of the reverse inequality is a consequence of Theorem \ref{CT0} and truly relies on the curvature assumption.
\par
A second consequence of Theorem \ref{CT0} is that on proper $\RCD(K,\infty)$ spaces the relaxed area functional coincides with the area functional itself. This is the content of Theorem \ref{CT2}. \par 
Let $\Omega \subset \X$ be an open bounded set and let $u \in \BV(\Omega)$. For every open set $E \subset \Omega$ we denote by $\Lip(E)$ the set of Lipschitz functions on $E$ and given $f \in \Lip(E)$ we denote by $\lip(f)$ its local Lipschitz constant. The relaxed area functional is defined by
\[
\tilde{\A}(u,E):=
\inf \Big\{ \liminf_{n \to + \infty} \int_E \sqrt{1+\lip(f_n)^2} \, d \m \Big| (f_n)_n \subset \sL^1(E) \cap \Lip(E), \|f_n-u\|_{\sL^1(E)} \to 0 \Big\}.
\]
\begin{customthm}{3} \label{CT2}
    Let $(\X,\sd,\m)$ be a proper $\RCD(K,\infty)$ space, let $\Omega \subset \X$ be open and bounded and let $u \in \BV(\Omega)$. Then, for every open set $E \subset \Omega$, it holds
    \[
    \tilde{\A}(u,E)=\A(u,E).
    \]
\end{customthm}

The previous theorem extends the analogous result proved in \cite{Tens} for Sobolev functions. Note that in general metric measure spaces (even if they are doubling and support a Poincaré inequality) the relaxed area functional does not coincide with the area functional itself, as shown in \cite{relax}. 
Theorem \ref{CT2}, in particular, allows to apply the mild regularity results for minimizers of the relaxed area functional in metric measure spaces established in \cite{HKLreg} to the area functional itself in the $\RCD(K,N)$ setting.
Next, we obtain a uniqueness result for perimeter minimizers (which follows through the identification between area and perimeter). In what follows, when we refer to an $\RCD(K,N)$ space, we implicitly mean that $N<+\infty$.

\begin{customthm}{4} \label{CT4}
    Let $(\X,\sd,\m)$ be a proper $\RCD(K,\infty)$ space, let $\Omega \subsetneq \X$ be an open set and let $u \in \sL^\infty(\Omega) \cap \W^{1,1}(\Omega)$ be a function whose hypograph is perimeter minimizing.
    If $v \in \sL^\infty(\Omega)$ has hypograph which is perimeter minimizing and $\{u \neq v\} \subset \subset \Omega$, then $u=v$.
\end{customthm}

\begin{corollary*}
    Let $(\X,\sd,\m)$ be an $\RCD(K,N)$ space, let $\Omega \subsetneq \X$ be an open set and let $u \in \W^{1,1}(\Omega)$ be a function whose hypograph is perimeter minimizing.
    If $v \in \sL^1(\Omega)$ has hypograph which is perimeter minimizing and $\{u \neq v\} \subset \subset \Omega$, then $u=v$.
\end{corollary*}

We then consider sequences of $\RCD(K,N)$ spaces and, relying on Theorems \ref{CT0} and \ref{CT4}, we show that Sobolev minimizers of the area functional in a limit space can be approximated with minimizers along the converging sequence of spaces.
This is the content of Theorem \ref{T|approssimazione aggiunta dopo} below. 
Theorem \ref{T|approssimazione aggiunta dopo} is useful when working with Ricci limit spaces, because it allows to pass information from minimizers of the area on smooth manifolds (for which there is a rich theory) to minimizers on the limit space. 
We work with functions that minimize the area (see Definition \ref{D|area minimizers}), but we could have worked with functions whose hypographs are perimeter minimizing since the two notions are equivalent (see Proposition \ref{P|equivalence}).

When considering pointed measured Gromov-Hausdorff convergence (pmGH) of spaces, we adopt the extrinsic viewpoint, meaning that we implicitly assume all spaces to be embedded in a common metric space $(Z,\sd_Z)$ where the convergence is realized. In particular, balls in this setting are to be intended w.r.t. the distance in $Z$ (where all spaces are embedded), while we will specify from time to time which measure is to be considered.
We refer to the beginning of Section \ref{S|approximation} for more details on the notation used.

\begin{customthm}{5} \label{T|approssimazione aggiunta dopo}
     Let $(\X_i,\sd_i,\m_i,x_i)$ be a sequence of $\RCD(K,N)$ spaces converging in pmGH sense to $(\X,\sd,\m,x)$. Let $r>0$ be such that $B_r(x) \subsetneq \X$ and let $u \in \W^{1,1}(B_r(x),\m)$ be an area minimizer. 
    For every $s \in (0,r)$, modulo passing to a subsequence, there exist area minimizers $u_i \in \BV(B_s(x),\m_i)$ such that $u_i 1_{B_s(x)}$ converge in $\sL_1$ to $u 1_{B_s(x)}$ and $\A(u_i,B_s(x),\m_i) \to \A(u,B_s(x),\m)$.
    \end{customthm}

Finally, in the Appendix, we collect some applications of the previously obtained theorems in the setting of Ricci limit spaces.

\medskip
\noindent
\textbf{Acknowledgements.}
I would like to thank Andrea Mondino for his guidance and Daniele Semola for useful comments on a previous version of the paper.

\section{Preliminaries} \label{S4}
We will work on metric measure spaces $\mms{X}{d}{m}$, where $(\ssf{X},\ssf{d})$ is a separable complete metric space and $\m$ is a non-negative Borel measure on $\ssf{X}$ which is finite on bounded sets. 
Given an open set $\Omega \subset \ssf{X}$ we denote by $\Lip(\Omega)$, $\Lip_c(\Omega)$ and $\Lip_b(\Omega)$ respectively the set of Lipschitz functions, Lipschitz functions with compact support and bounded Lipschitz functions on $\Omega$. Similarly, we denote by $C(\Omega)$ the set of continuous functions on $\Omega$ and by $\sL^0(\Omega)$ the one of measurable functions.
If $u \in \sL^0(\Omega)$ we also define
$\Hyp(u):=\{(x,t) \subset \Omega \times \bb{R}:t \leq u(x)\}$.
If $f \in \Lip(\Omega)$ and $x \in \Omega$ we define
\[
\lip(f)(x):= \limsup_{y \to x} \frac{|f(x)-f(y)|}{\ssf{d}(x,y)} \quad \text{and} \quad
\ssf{L}(f):= \sup_{x,y \in \Omega} \frac{|f(x)-f(y)|}{\ssf{d}(x,y)}.
\]
We say that $\Omega'$ is well contained in $\Omega$, and we denote it by $\Omega' \subset \subset \Omega$, if $\Omega'$ is bounded and $\sd(\Omega',{^c}\Omega)>0$. When we consider the spaces $\sL^1_{loc}(\Omega)$ and $\Lip_{loc}(\Omega)$ we implicitly mean that we are localizing w.r.t. well contained subsets of $\Omega$; note that these are not necessarily precompact in $\Omega$ (although they will be precompact if $\X$ is a finite dimensional $\RCD(K,N)$ space or a proper $\RCD(K,\infty)$ space). We now introduce Sobolev spaces in the setting of metric measure spaces, the main references being \cite{Cheeger}, \cite{AGS}, \cite{AGS2} and \cite{Giglimem}.

\begin{definition}
    Let $\mms{X}{d}{m}$ be a metric measure space and let $\Omega \subset \ssf{X}$ be an open set. A function $f \in \sL^2(\Omega)$ is said to be in the Sobolev space $\W^{1,2}(\Omega)$ if there exists a sequence $\{f_i\}_{i \in \bb{N}} \subset \Lip(\Omega)$ converging to $f$ in $\sL^2(\Omega)$ such that
    \[
    \limsup_{i \to + \infty} \int_\Omega \lip(f_i)^2 \, d \m < + \infty.
    \]
\end{definition}

For every $f \in \W^{1,2}(\Omega)$ one can define a function $|\nabla f| \in \sL^2(\Omega)$ such that for every open set $A \subset \Omega$ we have
\[
    \int_A |\nabla f|^2 \, d \m =
    \inf \Big\{ \liminf_{n \to + \infty} \int_A \lip(f_n)^2 \, d \m \Big| (f_n)_n \subset \sL^2(A) \cap \Lip(A), \|f_n-f\|_{\sL^2(A)} \to 0 \Big\}.
    \]
    The quantity in the previous expression will be called Cheeger energy and will be denoted by $\Ch(f)$ while $|\nabla f|$ will be called relaxed gradient.
    We define $\|f\|_{\W^{1,2}(\Omega)}:=\|f\|_{\sL^2(\Omega)}+\Ch(f)$. One can check that with this norm the space $\W^{1,2}(\Omega)$ is Banach.
    
\begin{remark} \label{R|asymptotic}
    Given a metric measure space $(\X,\sd,\m)$, an open set $\Omega \subset \X$ and a function $f \in \W^{1,2}(\Omega)$, there exists a sequence $f_n \in \Lip(\Omega)$ such that $\lip_a(f_n) \to |\nabla f|$ in $\sL^2(\Omega)$, where
    \[
    \lip_a (f)(x):= \lim_{r \to 0} \sup_{y,z \in B_r(x)} \frac{|f(y)-f(z)|}{\sd(y,z)}.
    \]
    This fact follows since in the definition of Sobolev spaces one could replace local Lipschitz constants $\lip$ with asymptotic Lipschitz constants $\lip_a$ and still obtain the same space (see \cite{ACDM15} and \cite[Remark $5.12$]{AGS}).
\end{remark}
    
We now introduce functions of bounded variation following \cite{Mir} (see also \cite{ADM}).

\begin{definition}
    Let $\mms{X}{d}{m}$ be a metric measure space and let $\Omega \subset \ssf{X}$ be an open set. A function $f \in \sL^1(\Omega)$ is said to be of bounded variation if there exists a sequence of $\{f_i\}_{i \in \bb{N}} \subset \Lip(\Omega)$ converging to $f$ in $\sL^1(\Omega)$ such that
    \[
    \limsup_{i \to + \infty} \int_\Omega \lip(f_i) \, d \m < + \infty.
    \]
    The space of such functions is denoted by $\BV(\Omega)$.
\end{definition}

    For every $f \in \BV(\Omega)$ and any open set $A \subset \Omega$ we define
    \[
    |Df|(A) =
    \inf \Big\{ \liminf_{n \to + \infty} \int_A \lip(f_n) \, d \m \Big| (f_n)_n \subset \sL^1(A) \cap \Lip(A), \|f_n-f\|_{\sL^1(A)} \to 0 \Big\}.
    \]
    One can check that the quantity in the previous expression is the restriction to the open subsets of $\Omega$ of a finite measure. We denote by $|\nabla f|$ the density of $|Df|$ with respect to $\m$ and by $|D^s f|$ the singular part of $|D u|$. We define $\|f\|_{\BV(\Omega)}:=\|f\|_{\sL^1(\Omega)}+|Df|(\Omega)$. One can check that with this norm the space $\BV(\Omega)$ is Banach. The function $f$ belongs to $\W^{1,1}(\Omega)$ if $f \in \BV(\Omega)$ and $|Df| << \m$.

\begin{definition}
    Let $(\X,\sd,\m)$ be a metric measure space and let $\Omega \subset \X$ be open. We define the perimeter of $E \subset \Omega$ in $\Omega$ by
    \[
    \sP(E,\Omega):=
    \inf \Big\{ \liminf_{n \to + \infty} \int_\Omega \lip(f_n) \, d \m \Big| (f_n)_n \subset \sL^1_{loc}(\Omega) \cap \Lip(\Omega), \|f_n-1_E\|_{\sL^1_{loc}(\Omega)} \to 0 \Big\}.
    \]
\end{definition}

We recall that $\sP(E,\cdot)$ is the restriction to open sets of a Borel measure and that if $E \subset \Omega$ has finite measure, then it has finite perimeter if and only if $1_E \in \BV(\Omega)$ and in this case $|D 1_E|=\sP(E,\cdot)$.

\begin{definition}
    We say that a metric measure space $\mms{X}{d}{m}$ is infinitesimally Hilbertian if the space $\W^{1,2}(\ssf{X})$ is a Hilbert space.
\end{definition}

We will consider pointed Gromov-Hausdorff (pGH) and pointed measured Gromov-Hausdorff (pmGH) convergence of metric measure spaces, and we refer to \cite{GMS13} for the relevant background.
We recall that in the case of a sequence of uniformly locally doubling metric measure spaces $(\ssf{X}_i,\ssf{d}_i, \m_i,\ssf{x}_i)$
(as in the case of $\RCD(K, N)$ spaces), pointed measured Gromov-Hausdorff convergence to $\pmms{X}{d}{m}{x}$ can be equivalently characterized asking for the existence of a proper
metric space $(Z, \ssf{d}_z )$ such that all the metric spaces $(\ssf{X}_i
, \ssf{d}_i)$ are isometrically embedded
into $(Z, \ssf{d}_z )$, $\ssf{x}_i \to \ssf{x}$ and $\m_i \to \m$ weakly in duality with continuous boundedly supported functions in $Z$. In this case we say that the convergence is realized in the space $Z$ (see \cite{GMS13}). 

\subsection{Properties of \texorpdfstring{$\RCD$}{RCD} spaces}

We now recall some properties of $\RCD(K,\infty)$ spaces; these were introduced in \cite{AGS2} (see
also \cite{AGMR,Giglimem}) coupling the Curvature Dimension condition $\ssf{CD}(K,\infty)$, previously pioneered
in \cite{S1,S2} and independently in \cite{V}, with the infinitesimal Hilbertianity assumption.
Whenever we consider an $\RCD(K,\infty)$ space $(X,\sd,\m)$, we assume that the support of $\m$ is the whole $X$.
The next proposition is an isoperimetric inequality that follows from the corresponding Poincaré inequality, which is proved in \cite{PI}.

\begin{proposition} 
Let $(\X,\sd,\m)$ be an $\RCD(K,\infty)$ space and let $E \subset B_{2r}(x)$ be a set of finite perimeter. Then we have the following isoperimetric inequality:
    \[
    \min \{\m(E \cap B_r(x)), \m(B_r(x) \setminus E) \} \leq 4r e^{|K|r^2} \sP(E,B_{2r}(x)).
    \]
\end{proposition}

Given $f :\X \times \bb{R} \to \bb{R}$ and $(x,t) \in \X \times \bb{R}$ we denote by $f^t$ and $f^x$ respectively the restriction of $f$ to $\X \times \{t\}$ and to $\{x\} \times \bb{R}$. The next proposition concerns a tensorization property of infinitesimally Hilbertian spaces and can be found in \cite{Tens} and \cite{AGS2}. We denote by $\sd_\times$ and $\m_\times$ respectively the product distance and the product measure in the space $\X \times \bb{R}$.

\begin{proposition} \label{P1}
Let $(\X,\sd,\m)$ be an $\RCD(K,\infty)$ space and let $f \in \Lip(\X \times \bb{R})$. Then we have
    \[
    |\nabla f|^2(x,t)=|\nabla f^x|^2(t)+|\nabla f^t|^2(x) \quad \text{for } \m_\times \text{-a.e. } (x,t) \in \X \times \bb{R}.
    \]
\end{proposition}

We now recall some properties of the heat flow in the $\RCD$ setting, referring to \cite{AGMR,AGS2} for the proofs of these results.
Given an $\RCD(K,\infty)$ space $(\X,\sd,\m)$, the heat flow $P_t:\sL^2(\X) \to \sL^2(\X)$ is the $\sL^2(\X)$-gradient flow of the Cheeger energy $\Ch$.
Let $\mathcal{P}_2(\X)$ denote the space of probability measures with finite second moment.
A key result in the $\RCD(K,\infty)$ theory is that the dual heat semigroup $\bar{P}_t: \mathcal{P}_2(\X) \to \mathcal{P}_2(\X)$ defined by
\[
\int_\X \, d \bar{P}_t \mu :=\int_\X P_tf \, d\mu \quad \forall \mu \in \mathcal{P}_2(\X), \quad \forall f \in \Lip_{b}(\X),
\]
is $K$-contractive with respect to the Wasserstein distance $W_2$ and for $t>0$ maps probability measures into absolutely continuous probability measures.
This allows to define a stochastically complete heat kernel $p_t:\X \times \X \to [0,+\infty)$, so that
the definition of $P_t(f)$ can be then extended to $\sL^{\infty}$ functions by setting
\[
P_t(f)(x):=\int_\X f(y) p_t(x,y) \, d\m(y).
\]
The heat flow has good approximation properties, in particular if $f \in \W^{1,2}(\X)$, then $P_t (f) \to f$ in $\W^{1,2}(\X)$; while if $f \in \sL^{\infty}(X)$, then $P_t f \in \Lip(\X)$ for every $t>0$.


In \cite{GBH14} it was proved that, for a \emph{proper} $\RCD(K,\infty)$ space $(\X,\sd,\m)$, if $f \in \BV(\X)$ then
\begin{equation} \label{E31}
|D P_t(f)| \leq e^{-Kt}\bar{P}_t(|Df|).
\end{equation}
Combining the previous contraction estimate with a lower semicontinuity argument we obtain that, if $\X$ is proper and
$\Omega \subset \X$ is open with $|Df|(\partial \Omega)=0$, then
    \[
    \lim_{t \to 0} |DP_t(f)|(\Omega)=|Df|(\Omega).
    \]

We next recall two basic lemmas on measures in metric measure spaces. Note in particular that these hold without the assumption that the space is proper. Given an $\bb{R}^m$-valued measure $\mu$ on $\X$ we denote by $|\mu|$ its variation. The first part of the next lemma follows from \cite[Proposition $1.23$]{AmbFuscPall}, while the second part follows repeating the argument of \cite[Proposition $1.47$]{AmbFuscPall} using the density of continuous functions in $\sL^1(\X)$ when $(\X,\sd,\m)$ is $\RCD(K,\infty)$.

\begin{lemma} \label{L2}
    Let $(\X,\sd,\m)$ be an $\RCD(K,\infty)$ space and let $f \in \sL^1(\X)^m$. Consider the $\bb{R}^m$-valued measure $\mu:=f\m$, then for every Borel set $B \subset \X$
    \[
    |\mu|(B)=\int_B|f| \, d\m,
    \]
    and for every open set $A \subset \X$ 
    \[
    |\mu|(A)=\sup_{\substack{\{a_i\}_{i=1}^m \subset C(A) \\ \sum a_i^2 \leq 1}} \sum_{i=1}^m \int_A a_i f_i \, d\m.
    \]
\end{lemma}

The next lemma follows by standard measure theoretic arguments.

\begin{lemma} \label{L1}
Let $(\X,\sd,\m)$ be an $\RCD(K,\infty)$ space
    and let $\mu,\nu$ be positive Radon measures on $\X$. If $\mu(A) \geq \nu(A)$ for every open set $A \subset \X$, then $\mu(E) \geq \nu(E)$ for every Borel set.
\end{lemma}

We now recall some properties of $\RCD(K,N)$ spaces.
The finite dimensional $\RCD(K,N)$ condition is obtained coupling the finite dimensional Curvature Dimension condition $\ssf{CD}(K,N)$ with the infinitesimal Hilbertianity assumption and was formalised in \cite{Giglimem}. For a thorough introduction to the topic we refer to the survey \cite{Asurv} and the references therein. 
Let us also mention that in the literature one can find also the (a-priori weaker) $\RCD^*(K,N)$. It was proved in \cite{EKS, AMS15}, that $\RCD^*(K,N)$ is equivalent to the dimensional Bochner inequality. Moreover, \cite{CavallettiMilman} (see also \cite{Liglob}) proved that $\RCD^*(K,N)$ and $\RCD(K,N)$ coincide. We now recall some results on $\RCD(K,N)$ spaces that we will need later in the paper.
\par 
The $\RCD(K,N)$ condition implies that the measure is locally doubling (see \cite{S1}) and, as mentioned before, the validity of a Poincaré inequality. In particular, if $f$ is a locally Lipschitz function on an $\RCD(K,N)$ space, its relaxed gradient coincides with its local Lipschitz constant $\lip(f)$ (see \cite[Theorem $12.5.1$]{HKST} after \cite{Cheeger}).

The next theorem is taken from \cite[Theorem $1.2$]{DephilGigli}. Given a metric space $(\X,\sd)$, we denote by $\aH^n$ the $n$-dimensional Hausdorff measure relative to $\sd$.

\begin{thm} \label{T|Gigli De philippis}
    Let $(\X_i,\sd_i,\aH^n_i,x_i)$ be a sequence of $\RCD(K,n)$ spaces such that $(\X_i,\sd_i,x_i)$ converges in pGH sense to a metric space $(\X,\sd,x)$. Then precisely one of the following holds:
    \begin{enumerate}
        \item 
        \[
        \limsup_{i \to + \infty} \aH^n_i(B_1(x_i))>0.
        \]
        In this case the $\limsup$ is a limit and the pGH convergence can be improved to pmGH convergence to $(\X,\sd,\aH^n,x)$. \\
        \item 
        \[
        \lim_{i \to + \infty} \aH^n_i(B_1(x_i))=0.
        \]
        In this case the Hausdorff dimension of $\X$ is bounded above by $n-1$.
    \end{enumerate}
\end{thm}

We now recall the definition of $\sL^1$-convergence of functions along a sequence of spaces converging in pmGH sense. 
As anticipated, we use the so-called extrinsic viewpoint and, when considering spaces $\X_i \to \X$ in pmGH sense, we work in the space $(Z,\sd_Z)$ which realizes the convergence. In particular, whenever we consider a ball $B_r(x)$ in this context, unless otherwise specified, we mean the ball of radius $r$ centered in $x$ in the space $Z$. 
Similarly, every function $u_i:\X_i \to \bb{R}$ will be considered, unless otherwise specified, as a function on $Z$ (extending it to zero outside of the embedding of $\X_i$ in $Z$). To avoid confusion we will specify from time to time which measure is to be considered when working with objects defined on $Z$.

\begin{definition} \label{D1}
    Let $(\X_i,\sd_i,\m_i,\ssf{x}_i)$ be a sequence of $\RCD(K,N)$ spaces converging in pmGH sense to $(\ssf{Y},\sd,\m,\ssf{y})$. We say that the functions $f_i \in \sL^1(\X_i,\m_i)$ converge to $f \in \sL^1(\X,\m)$ in $\sL^1$-sense if 
    \[
    \sigma \circ f_i \m_i \rightharpoonup \sigma \circ f\m 
    \quad 
    \text{and}
    \quad
    \int_{\X_i}|f_i| \, d\m_i \to \int_{\X}|f| \, d\m,
    \]
    where $\sigma(z):= \text{sign} \sqrt{|z|}$ and weak convergence is intended w.r.t. boundedly supported functions in $(Z,\sd_z)$, which is the space realizing the convergence.  We say that $f_i \in \sL^1_{loc}(\X_i,\m_i)$ converge in $\sL^1_{loc}$-sense to  $f \in \sL^1_{loc}(\X,\m)$ if $f_i 1_{B_r(x_i)} \to f1_{B_r(x)}$ in $\sL^1$-sense for every $r>0$.
\end{definition}

In the previous definition, $\sigma$ is needed to recover the usual $\sL^1$ convergence when we restrict to a fixed measured space $(X,\m)$. Indeed, since $\sL^1(X)$ in not uniformly convex in general, weak convergence together with convergence of the norms does not imply strong convergence. On the other hand, if $f_i \to f$ in $L^1$ sense on a fixed space $(X,\m)$ according to Definition \ref{D1}, then $\sigma \circ f_i \to \sigma \circ f$ in $\sL^2(X)$, and one can check that this implies $\int_X|f_i -f| \, d\m \to 0$.


The next proposition is taken from \cite[Proposition $3.3$]{ABS19}.
\begin{proposition} \label{P|compactness}
    Let $(\X_i,\sd_i,\m_i,\ssf{x}_i)$ be a sequence of $\RCD(K,N)$ spaces converging in pmGH sense to $(\ssf{Y},\sd,\m,\ssf{y})$. Let $f_i \subset \BV(\X_i,\m_i)$ with $\supp(f_i) \subset B_r(x) \subset Z$, where $(Z,\sd_Z)$ is the space realizing the convergence, be functions with
    \[
    \sup_{i \in \bb{N}} \Big( \|f_i\|_{\BV(\X_i,\m_i)}
    +
    \|f_i\|
    _{\sL^\infty(\X_i,\m_i)} \Big)< 
    + \infty.
    \]
    Then there exists a (non relabeled) subsequence and $f \in \BV(\X,\m)$ with $\supp{f} \subset \bar{B}_r(x)$ such that $f_i \to f$ in $\sL^1$.
\end{proposition}

The next proposition is taken from \cite[Proposition $3.6$]{ABS19}.

\begin{proposition} \label{P|lsc variation aperti}
    Let $(\X_i,\sd_i,\m_i,x_i)$ be a sequence of $\RCD(K,N)$ spaces converging in pmGH sense to $(\X,\sd,\m,x)$. 
    If $u \in \BV(\X,\m)$, and $u_i \in \BV(\X_i,\m_i)$ is a sequence such that $u_i \to u$ in $\sL^1$ 
    and $\sup_i \|u_i\|_{\sL^\infty(\X_i,\m_i)} < + \infty$,
    then for every open set $A \subset Z$, where $(Z,\sd_Z)$ is the metric space realising the convergence, we have
    \[
    |Du|(A) \leq \liminf_{i} |Du_i|(A).
    \]
\end{proposition}

\subsection{Minimal sets}

We now turn our attention to minimal sets. We only recall the properties that will be used later in the paper.

\begin{definition}
Let $(\X,\sd,\m)$ be a metric measure space and let $\Omega \subset \X$ be an open set.
    Let $E \subset \Omega$ be a set of locally finite perimeter. We say that $E$ is perimeter minimizing in $\Omega$ if for every $x \in \Omega$, $r>0$ and $F \subset  \Omega$ such that $F \Delta E \subset \subset B_r(x) \cap \Omega$ we have that $P(E,B_r(x) \cap \Omega) \leq P(F,B_r(x) \cap \Omega)$. If we say that $E$ is perimeter minimizing we implicitly mean that $\Omega=\X$.
\end{definition}


The next theorem comes from \cite[Theorem $4.2$ and Lemma $5.1$]{Dens}.

\begin{thm} \label{T|density}
    Let $(\X,\sd,\m)$ be an $\RCD(K,N)$ space and let $\Omega \subset \X$ be an open set.
    There exist $C,\gamma>0$ depending only on $K$ and $N$ such that the following hold. If $E \subset \X$ is a set minimizing the perimeter in $\Omega \subset \X$, then, up to modifying $E$ on an $\m$-negligible set, for any $x \in \partial E$ and $r>0$ such that $B_{2r}(x) \subset \Omega$ we have
    \[
    \frac{\m(E \cap B_r(x))}{\m(B_r(x))}>\gamma, \quad \frac{\m(B_r(x) \setminus E)}{\m(B_r(x))}>\gamma
    \]
    and
    \[
    \frac{\m(B_r(x))}{Cr} \leq P(E,B_r(x)) \leq \frac{C\m(B_r(x))}{r}.
    \]
\end{thm}


The next proposition is taken from \cite[Theorem $2.43$]{Weak}.

\begin{proposition} \label{P26}
Let $(\X_i,\sd_i,\m_i,\ssf{x}_i)$ be a sequence of $\RCD(K,N)$ spaces converging in pmGH sense to $(\ssf{Y},\sd,\m,\ssf{y})$.
    Let $E_i \subset \X_i$ be a sequence of Borel sets converging in $\sL^1_{loc}$-sense to $E \subset \ssf{Y}$.
    Assume that each $E_i$ is perimeter minimizing in $B_{r_i}(\ssf{x}_i)$ and that $r_i \uparrow + \infty$. Then $E$ is perimeter minimizing and in the metric space realizing the convergence we have that $\partial E_i \to \partial F$ in Kuratowski sense.
\end{proposition}

We recall that given a space $(\X,\sd,\m)$
we denote by $\sd_\times$ and $\m_\times$ respectively the product distance and the product measure on $\X \times \bb{R}$.
Moreover, given a function $u \in \Lip_{loc}(\X)$ we will use the following notation for balls on the graph of $u$ (these are actually intersections of balls in the product with the graph of $u$):
\[
B^m_r(x,t):=B^\times_r(x,t) \cap \graph(u) \subset \X \times \bb{R}.
\]
Similarly, given a point $x \in \X$, we will denote by $\bar{x}$ its projection on the graph of $u$, i.e. $\bar{x}:=(x,u(x))$.
Whenever we evaluate a function $f:\X \to \bb{R}$ on the graph of $u$ we implicitly mean that we are referring to 
$f \circ \pi$, where $\pi:\graph(u) \to \X$ is the standard projection.

The next result follows from \cite[Theorem $5.8$]{C2}.

\begin{thm} \label{T|Harnack}
    There exists a constant $C>0$ depending only on $N$ such that if
      $(\X,\sd,\m)$ is an $\RCD(0,N)$ space and $u \in \Lip_{loc}(B_{3R}(p))$ is a positive function whose hypograph minimizes the perimeter in $B_{3R}(p) \times \bb{R}$, setting $\bar{p}:=(p,u(p))$,  we have
     \[
         \sup_{B^m_{R/2}(\bar{p})} u \leq C \inf_{ B^m_{R/2}(\bar{p})} u.
    \]
\end{thm}

\section{Heat flow approximation of the area and applications in \texorpdfstring{$\RCD(K,\infty)$}{RCD(K,8)} spaces}

In this section, we prove Theorems \ref{CT0} to \ref{CT4}. We denote by $\sd_\times$ and $\m_\times$ respectively the product distance and the product measure on $\X \times \bb{R}$.
The next two propositions contain the inequality in the area formula that was proved in \cite{Tens} for $PI$ spaces. Similar proofs work in our setting and we repeat them (with the few due modifications) for the sake of completeness as they will be needed to prove Theorems \ref{CT0} and \ref{CT1}.


\begin{definition} \label{D|Area functional}
    Let $(\X,\sd,\m)$ be an $\RCD(K,\infty)$ space, let $\Omega \subset \X$ be an open set and let $u \in \BV(\Omega)$. For every measurable $E \subset \Omega$ we define the area of $u$ on $E$ as
    \[
    \A(u,E):=\int_E \sqrt{1+|\nabla u|^2} \, d\m + |D^s u|(E).
    \]
\end{definition}

\begin{proposition} \label{P2}
Let $(\X,\sd,\m)$ be an $\RCD(K,\infty)$ space and let $\Omega \subset \X$ be an open set. If $u \in \BV(\Omega)$, then for every Borel set $E \subset \Omega$, we have
    \[
\sP(\Hyp(u),E \times \bb{R}) \geq \A(u,E).
\]
\begin{proof}
    We will consider the measure $\mu:=P(\Hyp(u),\cdot \times \bb{R})$ on the Borel subsets of $\Omega$ and we will show separately that 
    \begin{equation} \label{E2}
    \mu \geq \sqrt{1+|\nabla u|^2} \m
\quad \text{and} \quad
    \mu \geq |D^s u|.
    \end{equation}
    Since $\m \perp |D^s u|$ this will then imply that $\mu \geq \sqrt{1+|\nabla u|^2} \m+|D^s u|$. \par 
    Thanks to Lemma \ref{L1}, to obtain \eqref{E2}, it is sufficient to show that for every open set $E \subset \Omega$ we have
    \[
    \mu(E) \geq \int_E \sqrt{1+|\nabla u|^2} \, d\m
\quad \text{and} \quad
    \mu(E) \geq |D^s u|(E).
    \]
    \par 
    We now repeat an argument from \cite{Tens} to obtain that, for every open set $E \subset \X$ and $a,b \in C(E)$ such that $a^2+b^2 \leq 1$, it holds
    \begin{equation} \label{E1}
    \mu(E) \geq \int_E a \, d\m + \int_E b \, d|Du|.
    \end{equation}
    We report the argument for the sake of completeness. Fix $N \in \bb{N}$ and let $g_n \in \Lip(E \times (-N,N))$ be convergent to $1_{\Hyp(u)}$ in $\sL^1(E \times (-N,N))$. Modulo passing to a subsequence, we may assume that $g_n^t \to 1_{\Hyp(u)}^t$ in $\sL^1(E)$ for $\lambda^1$-a.e. $t \in (-N,N)$ and that $g_n^x \to 1_{\Hyp(u)}^x$ in $\sL^1((-N,N))$ for $\m$-a.e. $x \in E$. For every $a,b \in C(E)$ with $a^2+b^2 \leq 1$, applying Proposition \ref{P1}, we get
    \begin{align*}
    \liminf_{n \to + \infty} & \int_{E \times (-N,N)} |\nabla g_n| \, d\m_\times \geq
    \liminf_{n \to + \infty} \int_{E \times (-N,N)} a(x)|\nabla g_n^x|(t)+b(x)|\nabla g_n^t|(x) \, d\m_\times (x,t) 
    \\
    & \geq \int_E \liminf_{n \to + \infty} \int_{(-N,N)}a(x)|\nabla g_n^x|(t) \, dt d\m(x)+
    \int_{(-N,N)} \liminf_{n \to + \infty} \int_\X b(x)|\nabla g_n^t|(x) \, d\m(x) dt
    \\
    & \geq \int_E a 1_{\{u \in (-N,N)\}} \, d\m+\int_{(-N,N)}\int_Eb \, d |D1_{\{u>t\}}| dt.
    \end{align*}
    Since $g_n$ is arbitrary, we get
    \[
    \sP(\Hyp(u),E \times (-N,N)) \geq \int_E a 1_{\{u \in (-N,N)\}} \, d\m+\int_{(-N,N)}\int_Eb \, d |D1_{\{u>t\}}| dt,
    \]
    so that letting $N \uparrow + \infty$ and using the coarea formula we deduce \eqref{E1}. \par
    Equation \eqref{E1} then implies that
    \[
    \mu(E) \geq \int_E a +  b |\nabla u| \, d\m.
    \]
    Taking the supremum over all couples $a,b$ of the aforementioned type and using Lemma \ref{L2}, we get
    \[
    \mu(E) \geq \int_E \sqrt{1+|\nabla u|^2} \, d\m.
    \]
    Similarly, taking $a \equiv 0$ and $b\equiv 1$ in \eqref{E1}, we obtain
    \[
    \mu(E) \geq |D u|(E) \geq |D^s u|(E),
    \]
    concluding the proof.
\end{proof}
\end{proposition}

\begin{proposition} \label{P4}
Let $(\X,\sd,\m)$ be an $\RCD(K,\infty)$ space and let $\Omega \subset \X$ be an open set. If $u \in \Lip(\Omega)$, then, for every Borel set $E \subset \Omega$, it holds
    \[
\sP(\Hyp(u),E \times \bb{R}) = \A(u,E).
\]
\begin{proof}
    Thanks to the previous proposition, we only need to show that
    \[
    \sP(\Hyp(u),E \times \bb{R}) \leq \int_E \sqrt{1+|\nabla u|^2} \, d\m.
    \]
    Thanks to Lemma \ref{L1}, we only need to verify the previous inequality under the assumption that $E$ is open. Let $\chi_\epsilon \in C^{\infty}(\bb{R})$ be monotonically convergent to $1_{(0,+\infty)}$, with values in $[0,1]$ and $\int_\bb{R} \chi_\epsilon' \, dt \leq 1$. Then the functions
    \[
    g_\epsilon(x,t):=\chi_\epsilon(u(x)-t)
    \]
    are Lipschitz in $\Omega \times \bb{R}$ for every $\epsilon>0$ and converge in $\sL^1_{loc}(\Omega \times \bb{R})$ to the characteristic function of $\Hyp(u)$. Moreover Proposition \ref{P1} gives
    \[
    |\nabla g_\epsilon|(x,t)=\chi_\epsilon'(u(x)-t)\sqrt{1+|\nabla u|^2(x)} \quad \text{for } \m_\times \text{-a.e. } (x,t) \in \X \times \bb{R}.
    \]
    Integrating both sides and using Fubini's Theorem, we get
    \[
    \int_{E \times \bb{R}} |\nabla g_\epsilon| \, d \m_\times \leq \int_E \sqrt{1+|\nabla u|^2} \, d\m,
    \]
    which implies 
    \[
    \sP(\Hyp(u),E \times \bb{R}) \leq \int_E \sqrt{1+|\nabla u|^2} \, d\m.
    \]
\end{proof}
\end{proposition}

We now prove Theorem \ref{CT0} from the Introduction. Here we add the properness assumption on the space  $(\X,\sd,\m)$ as we will need to use \eqref{E31}.

\begin{customthm}{1}
    Let $(\X,\sd,\m)$ be a proper $\RCD(K,\infty)$ space and let $u \in \sL^{\infty}(\X) \cap \BV(\X)$. Then, for every bounded open set $E \subset \X$ such that $|D^su|(\partial E)=0$,
it holds
\[
\lim_{t \to 0} \A(P_t(u),E) = \A(u,E).
\]
\begin{proof}
    For every $t > 0$, we set $u_t:=P_t(u)$ and we note that this function is Lipschitz. 
We first show that
\[
\int_{E} \sqrt{1+|\nabla u|^2} \, d\m +|D^su|(E)
\leq \liminf_{t \to 0} \int_{E} \sqrt{1+|\nabla u_t|^2} \, d\m.
\]
 Using Propositions \ref{P2} and \ref{P4} and the lower semicontinuity of perimeters we get
\begin{align*}
\int_{E} & \sqrt{1+|\nabla u|^2} \, d\m +|D^su|(E) \leq \sP(\Hyp(u),E \times \bb{R}) \\
 & \leq \liminf_{t \to 0} \sP(\Hyp(u_t),E \times \bb{R})=
\liminf_{t \to 0} \int_{E} \sqrt{1+|\nabla u_t|^2} \, d\m,
\end{align*}
as desired. \par 
We now show that 
\[
\int_{E} \sqrt{1+|\nabla u|^2} \, d\m +|D^su|(E)
\geq \limsup_{t \to 0} \int_{E} \sqrt{1+|\nabla u_t|^2} \, d\m.
\]
Note that
\[
\int_E \sqrt{1+|\nabla u_t|^2} \, d\m=
\sup_{\substack{(a,b) \in C(E) \times C(E) \\ a^2+b^2 \leq 1}}
\int_E a+|\nabla u_t| b \, d\m.
\]
For every such pair $(a,b)$, we have 
\begin{align*}
\int_E & a+|\nabla u_t| b \, d\m \leq 
\int_E a \, d\m + \int_E be^{-Kt} \, d\bar{P}_t(|Du|)
\\
& \leq \int_E a \, d\m + \int_E be^{-Kt} \, P_t(|\nabla u|) \, d\m+e^{-Kt}\bar{P}_t(|D^su|)(E).
\end{align*}
Hence, passing to the supremum with respect to the pairs $(a,b)$, we get
\[
\int_E \sqrt{1+|\nabla u_t|^2} \, d\m \leq 
\int_E \sqrt{1+e^{-2Kt}P_t(|\nabla u|)^2} \, d\m +
e^{-Kt}\bar{P}_t(|D^su|)(E).
\]
Moreover, by definition of $\bar{P}_t$, we have that $\bar{P}_t(|D^su|)$ converges weakly to $|D^su|$ preserving the mass for every $t$, and since $|D^su|(\partial E)=0$ this implies that $\bar{P}_t(|D^su|)(E) \to |D^su|(E)$. In particular, passing to the limit as $t \to 0$ in the previous inequality, we get
\[
\limsup_{t \to 0}
\int_E \sqrt{1+|\nabla u_t|^2} \, d\m \leq 
\int_E \sqrt{1+|\nabla u|^2} \, d\m +|D^su|(E).
\]
\end{proof}
\end{customthm}

The next two results are needed to prove Theorem \ref{CT1} from the Introduction.

\begin{thm} \label{T1}
Let $(\X,\sd,\m)$ be a proper $\RCD(K,\infty)$ space.
    Let $\Omega \subset \X$ be an open bounded set and let $u \in \BV(\Omega)$. For every Borel set $E \subset \Omega$ it holds
\[
\sP(\Hyp(u),E \times \bb{R})=\A(u,E).
\]
\begin{proof}
The inequality
\[
\sP(\Hyp(u),E \times \bb{R}) \geq \int_E \sqrt{1+|\nabla u|^2} \, d\m +|D^su|(E)
\]
is true by Proposition \ref{P2} so we only show the reverse one. \par 
    Suppose first that $u \in \sL^{\infty}(\Omega) \cap \BV(\Omega)$.
    It is clear that it is sufficient to consider the case when $E$ is well contained in $\Omega$ and such that $|Du|(\partial E)=0$ and then argue by approximation.
So we assume that $E$ has these properties and, modulo multiplying $u$ by a cut off function which is equal to $1$ in a neighborhood of $\bar{E}$, we can also suppose that $u \in \sL^{\infty}(\Omega) \cap \BV(\Omega)$ with $\supp(u) \subset \subset \Omega$.
In particular, when we extend $u$ to be $0$ outside of $\Omega$, we obtain a function in $u \in \sL^{\infty}(\X) \cap \BV(\X)$. When we refer to $u$ in the next lines we implicitly mean its extension to $\X$.
By lower semicontinuity of perimeters and the previous proposition we then get
\begin{align*}
\sP &(\Hyp(u),E \times \bb{R}) \leq \liminf_{t \to 0} \sP(\Hyp(P_t(u)),E \times \bb{R})
\\
& =
\liminf_{t \to 0}
\int_E \sqrt{1+|\nabla P_t(u)|^2} \, d\m =
\int_E \sqrt{1+|\nabla u|^2} \, d\m +|D^su|(E),
\end{align*}
concluding the proof when $u \in \sL^{\infty}(\Omega) \cap \BV(\Omega)$. 
\par
If $u \in \BV(\Omega)$ is not bounded, we consider for every $k>0$ the function $u^k:=-k \vee u \wedge k$ and we exploit the lower semicontinuity of perimeters to obtain that for every open set $E \subset \Omega$ we have
\begin{align*}
\sP &(\Hyp(u),E \times \bb{R}) \leq
\liminf_{k \to + \infty} \sP(\Hyp(u^k),E \times \bb{R})
\\
& = \liminf_{k \to + \infty} 
\int_E \sqrt{1+|\nabla u^k|^2} \, d\m +|D^su^k|(E).
\end{align*}
Moreover, since $|D u^k| \leq |Du|$ as measures, the last term is controlled by
\[
\int_E \sqrt{1+|\nabla u|^2} \, d\m +|D^su|(E),
\]
concluding the proof.
\end{proof}
\end{thm}

\begin{proposition} \label{P3}
    Let $(\X,\sd,\m)$ be a proper $\RCD(K,\infty)$ space. If $u : B_{2r}(x) \to \bb{R}$ is a measurable function whose hypograph has finite perimeter, then $u \in \BV(B_r(x))$. Similarly if $\Omega$ is bounded and $u \in \sL^1(\Omega)$ has hypograph of finite perimeter, then $u \in \BV(\Omega)$.
    \begin{proof}
    Assume first that $\Omega$ is bounded and that $u \in \sL^1(\Omega)$ is a function whose hypograph has finite perimeter. 
    Fix $N \in \bb{N}$ and let $g_n \in \Lip(\Omega \times (-N,N))$ be convergent to $1_{\Hyp(u)}$ in $\sL^1(\Omega \times (-N,N))$. Modulo passing to a subsequence, we may assume that $g_n^t \to 1_{\Hyp(u)}^t$ in $\sL^1(\Omega)$ for $\lambda^1$-a.e. $t \in (-N,N)$. Then
    \[
    \liminf_{n \to + \infty} \int_{\Omega \times (-N,N)} |\nabla g_n| \, d\m_\times \geq
    \liminf_{n \to + \infty} \int_{\Omega \times (-N,N)} |\nabla g_n^t| \, d\m_\times
    \geq \int_{-N}^N \sP(\{u>t\},\Omega) dt.
    \]
    Since $g_n$ is arbitrary, letting $N \uparrow + \infty$, we deduce 
    \begin{equation} \label{Eff}
    \int_{\bb{R}} \sP(\{u>t\},\Omega) dt < + \infty.
    \end{equation}
    Since $u \in \sL^1(\Omega)$, this is sufficient to conclude that $u \in \BV(\Omega)$. \par
    If $u :B_{2r}(x) \to \bb{R}$ is measurable
    we need to prove that $u \in \sL^1(B_r(x))$, and then we get that $u \in \BV(B_r(x))$ by the previous part of the proof. To this aim we write
    \begin{align*}
    \int_{B_r(x)}|u| & \, d\m=
    \int_0^{\infty} \m(\{|u| \geq t\} \cap B_r(x)) \, dt
    \\
    & \leq T\m(B_r(x))+\int_T^{\infty} \m(\{|u| \geq t \cap B_r(x)\}) \, dt.
    \end{align*}
    If $T \in \bb{N}$ is large enough, for every $t \geq T$ the isoperimetric inequality gives 
    \[
    \m(\{|u| \geq t \cap B_r(x)\}) \leq c(K,r)\sP(\{|u| \geq t\},B_{2r}(x)).
    \]
    In particular, combining this with \eqref{Eff}, we get
    \begin{align*}
    \int_{B_r(x)}|u| \, d\m
    \leq
    T\m(B_r(x))+
    c(K,r)\int_T^{\infty} \sP(\{|u| \geq t\},B_{2r}(x)) \, dt
    < + \infty,
    \end{align*}
    concluding the proof.
    \end{proof}
\end{proposition}

Combining Theorem \ref{T1} and the previous proposition we immediately get Theorem \ref{CT1}, which we recall below.

\begin{customthm}{2} 
Let $(\X,\sd,\m)$ be a proper $\RCD(K,\infty)$ space, let $\Omega \subset \X$ be open and bounded and let $u \in \sL^1(\Omega)$. Then $\sP(\Hyp(u),\Omega \times \bb{R})<+\infty$ if and only if $u \in \BV(\Omega)$. In this case, for every Borel set $E \subset \Omega$, it holds
    \[
\sP(\Hyp(u),E \times \bb{R})=\A(u,E).
\]
Moreover, if $u \in \sL^0(B_{2r}(x))$ and $\sP(\Hyp(u),B_{2r}(x) \times \bb{R})<+\infty$, then $u \in \BV(B_r(x))$ and the previous equation holds for every Borel set $E \subset B_r(x)$.
\end{customthm}

We now turn our attention to Theorem \ref{CT2}, which will be an easy consequence of the approximation property of Theorem \ref{CT0}.
As anticipated in the Introduction, given $u \in \BV(\Omega)$ and an open set $E \subset \Omega$, we define the relaxed area functional as
\[
\tilde{\A}(u,E):=
\inf \Big\{ \liminf_{n \to + \infty} \int_E \sqrt{1+\lip(f_n)^2} \, d \m \Big| (f_n)_n \subset \sL^1(E) \cap \Lip(E), \|f_n-u\|_{\sL^1(E)} \to 0 \Big\}.
\]
In \cite{Tens} it was proved that for every $f \in \W^{1,1}(\Omega)$, we have 
    \[
    \tilde{\A}(f,E)=\int_E \sqrt{1+|\nabla f|^2} \, d\m,
    \]
and in Theorem \ref{CT2} we extend this to functions in $\BV(\Omega)$.

\begin{customthm}{3}
     Let $(\X,\sd,\m)$ be a proper $\RCD(K,\infty)$ space, let $\Omega \subset \X$ be open and bounded and let $u \in \BV(\Omega)$. Then, for every open set $E \subset \Omega$, it holds
    \[
    \tilde{\A}(u,E)=\A(u,E).
    \]
    \begin{proof}
        We first show that
        for every open set $E \subset \Omega$ we have
        \[
        \tilde{\A}(u,E) \geq \int_E \sqrt{1+|\nabla u|^2} \, d\m +|D^su|(E).
        \]
        To this aim fix an open set $E$ and let $\{u_n\}_{n \in \bb{N}} \subset \Lip(E)$ be such that 
        \[
        \int_E \sqrt{1+|\nabla u_n|^2} \, d\m \to \tilde{\A}(u,E)
        \]
        and $u_n \to u$ in $\sL^1(E)$. By the lower semicontinuity of perimeters, we get
        \[
        \A(u,E)=
       \int_E \sqrt{1+|\nabla u|^2} \, d\m +|D^su|(E)  \leq \liminf_{n \to + \infty} \int_E \sqrt{1+|\nabla u_n|^2} \, d\m =\tilde{\A}(u,E).
        \]
        We now show that 
        \[
        \tilde{\A}(u,E) \leq \A(u,E)
        \]
        for every open set $E \subset \Omega$.
        To this aim it is sufficient to consider open sets $E \subset \subset \Omega$ such that $|D^s u|(\partial E)=0$ and then argue by approximation. Let $u^k:=-k \vee u \wedge k$ and note that for every $\epsilon>0$ we have that if $k$ is sufficiently large
        $
        \tilde{\A}(u)(E) \leq \tilde{\A}(u^k)(E)+\epsilon.
        $
        Using the lower semicontinuity of $\tilde{\A}$, the fact that $\tilde{\A}$ is the area functional on Lipschitz functions and Theorem \ref{CT0}, we then obtain that
        \begin{align*}
         \tilde{\A}(u^k,E) \leq 
         \liminf_{t \to 0} 
         \tilde{\A}(P_t(u^k),E)
         & =
         \liminf_{t \to 0}\int_E \sqrt{1+|\nabla P_t(u^k)|^2} \, d\m 
         \\
         & =\int_E \sqrt{1+|\nabla u^k|^2} \, d\m +|D^su^k|(E).
        \end{align*}
        In particular we get that, if $k$ is large enough, we have
        \[
        \tilde{\A}(u,E) \leq\int_E \sqrt{1+|\nabla u^k|^2} \, d\m +|D^su^k|(E)+\epsilon,
        \]
        so that $\tilde{\A}(u,E) \leq \int_E \sqrt{1+|\nabla u|^2} \, d\m +|D^su|(E)+ \epsilon$. Since $\epsilon>0$ is arbitrarily small, we conclude.
    \end{proof}
\end{customthm}

We conclude the section with some basic properties of minimizers of the area functional in $\RCD(K,\infty)$ spaces that are needed to prove Theorem \ref{CT4} from the Introduction.

\begin{definition} \label{D|area minimizers}
    Let $(\X,\sd,\m)$ be an $\RCD(K,\infty)$ space, let $\Omega \subset \X$ be an open set and let $u \in \BV_{loc}(\Omega)$. We say that that $u$ is an area minimizer in $\Omega' \subset \Omega$ if for every $f \in \BV(\Omega')$ such that $\{f \neq u\} \subset \subset \Omega'$ we have
    \[
    \A(u,\Omega') \leq \A(f,\Omega').
    \]
    If we say that $u \in \BV_{loc}(\Omega)$ is an area minimizer, we mean that the previous condition is satisfied for every $\Omega' \subset \Omega$.
\end{definition}

The following proposition follows from \cite[Theorem $4.2$]{HKLreg} keeping in mind that in our setting the area functional coincides with the relaxed area functional thanks to Theorem \ref{CT2}. Since the setting of \cite{HKLreg} is the one of $PI$ spaces (and the proof of \cite[Theorem $4.2$]{HKLreg} relies crucially on the doubling assumption), the next result holds in the finite dimensional setting of $\RCD(K,N)$ spaces.

\begin{proposition} \label{P|local boundedness}
    Let $(\X,\sd,\m)$ be an $\RCD(K,N)$ space, let $\Omega \subset \X$ be an open set and let $u \in \BV_{loc}(\Omega)$ be an area minimizer, then $u \in \sL^\infty_{loc}(\Omega)$.
\end{proposition}

A similar result holds for functions whose hypograph is perimeter minimizing
adapting an argument that can be found in \cite{Giusti}. We report such argument with the due adaptations.

\begin{proposition} \label{P|aggiunta}
    Let $(\X,\sd,\m)$ be an $\RCD(K,N)$ space, let $\Omega \subset \X$ be open and let $u \in \sL^1(\Omega)$ be a function whose hypograph minimizes the perimeter. Then $u \in \sL^\infty_{loc}(\Omega)$. 
    \begin{proof} 
        Let $x \in \Omega$ and suppose for simplicity that $B_3(x) \subset \subset \Omega$. 
        We will show that $u$ is uniformly bounded in its Lebesgue points in $B_1(x)$, which will then imply the claim. To this aim let $y \in B_1(x)$ be a Lebesgue point for $u$ and assume without loss of generality that $u(y)>0$. 
        Denoting by $\Hyp(u)$ the closed representative of the hypograph, we have that $(y,u(y)) \in \partial \Hyp(u)$, so that the density estimates from Theorem \ref{T|density} imply that
        \[
        \m_\times (B_1^\times (y,u(y)) \cap \Hyp(u)) \geq c(K,N).
        \]
        In particular if $T \in \bb{N}$ and $u(y) > 2T$ we have
        \begin{align*}
        \m_\times(\Hyp(u) \cap B_1(y) \times [0,2T]) &\geq \sum_{i=1}^T  \m_\times(\Hyp(u) \cap B^\times_1(y,2i) ) \\
        & \geq T\m_\times (B_1^\times (y,u(y)) \cap \Hyp(u)) \geq T c(K,N).
        \end{align*}
        This implies that
        \[
        \int_{B_2(x)} u \vee 0 \, d\m \geq \int_{B_1(y)} u \vee 0 \, d\m \geq T c(K,N).
        \]
        Putting everything together, if $y \in B_1(x)$ is a Lebesgue point of $u$ and $u(y)>2T$, then $T \leq C(x,u,K,N) < + \infty$. This shows that for every such $y$ we have $u(y) \leq 2C(x,u,K,N)+1$, proving the claim and the statement.
    \end{proof}
\end{proposition}

The proof of the next proposition follows by adapting arguments from \cite[Section $3.2$]{C2} and we report it for the sake of completeness.

\begin{proposition} \label{P|equivalence}
    Given a proper $\RCD(K,\infty)$ space $(\X,\sd,\m)$, an open set $\Omega \subset \X$ and a function $u \in \sL^\infty(\Omega) \cap \BV(\Omega)$, then $u$ is a minimizer of the area functional in $\Omega$ if and only if its hypograph is a perimeter minimizer in $\Omega \times \bb{R}$. 
    \begin{proof}
          If the hypograph of $u$ is perimeter minimizing, then 
          for every $f \in \BV(\Omega) \cap \sL^\infty(\Omega)$ such that $\{f \neq u\} \subset \subset \Omega$ we have $\Hyp(u) \Delta \Hyp(f) \subset \subset \Omega \times \bb{R}$, so that applying Theorem \ref{CT1} we obtain
    \[
    \A(u,\Omega) \leq \A(f,\Omega).
    \]
    Similarly, if  $f \in \BV(\Omega) $ is such that $\{f \neq u\} \subset \subset \Omega$, then considering $f^c:=-c \vee f \wedge c$ with $c:=\|u\|_{\sL^\infty}+1$ we get by the previous case that
    \[
    \A(u,\Omega) \leq \A(f^c,\Omega) \leq \A(f,\Omega).
    \]
    We now prove the reverse implication, i.e. that if
    $u \in \sL^\infty(\Omega) \cap \BV(\Omega)$ is a minimizer of the area functional in $\Omega$, then its hypograph is a perimeter minimizer in $\Omega \times \bb{R}$. \par 
    To this aim, consider a set $E \subset \Omega \times \bb{R}$ such that $\Hyp(u) \Delta E \subset \subset \Omega \times \bb{R}$. Modulo translating vertically, we may suppose that there exists $c>1$ such that $u$ takes values in $(1,c-1)$ and $\Hyp(u) \Delta E \subset \subset \Omega \times (1,c-1)$.
    We then define $w(E):\Omega \to \bb{R}$ by
    \[
    w(E)(x):=\int_0^c 1_E(x,s) \, ds
    \]
    and we claim that
    \begin{equation} \label{E|w(E)}
    w(E) \in \BV(\Omega)
    \quad 
    \text{and}
    \quad
    \A(w(E),\Omega) \leq \sP(E,\Omega \times \bb{R}).
    \end{equation}
    If the claim holds, noting that by construction $w(E)$ is a competitor for $u$, we get
    \[
    \sP(\Hyp(u),\Omega \times \bb{R})=\A(u,\Omega) \leq \A(w(E),\Omega) \leq \sP(E,\Omega \times \bb{R}).
    \]
    To prove the claim, we consider a sequence $f_n \in \Lip(\Omega \times (0,c))$ converging in $\sL^1(\Omega \times (0,c))$ to $1_E$ and such that
    \[
    \lim_{n \to + \infty} |D f_n|(\Omega \times (0,c))=\sP(E,\Omega \times \bb{R}).
    \]
    Modulo truncating, we can assume that $f_n \equiv 1$ on $\Omega \times \{0\}$ and $f_n \equiv 0$ on $\Omega \times \{c\}$.
    We then define $w(f_n): \Omega \to \bb{R}$ as
    \[
    w(f_n)(x):=\int_0^c f_n(x,s) \, ds.
    \]
    These functions are Lipschitz since
    \[
    \frac{|w(f_n)(x)-w(f_n)(y)|}{\sd(x,y)} \leq \int_0^c \frac{|f_n(x,s)-f_n(y,s)|}{\sd(x,y)} \, ds \leq c \ssf{L}(f_n).
    \]
    We now use the notation defined before Proposition \ref{P1}.
    By reverse Fatou Lemma and the fact that each $f_n$ is Lipschitz, for $\m$-a.e. $x \in X$ it holds
    \begin{equation} \label{E30}
    |\nabla w(f_n)|(x)=\lip (w(f_n))(x) \leq \int_0^c |\nabla f_n^s|(x) \, ds.
    \end{equation}
    Moreover, using the tensorization property of Proposition \ref{P1}, it holds
    \begin{align*}
    |D f_n|(\Omega \times (0,c))
    & =
    \int_{\Omega \times (0,c)} |\nabla f_n |(x,s) \, d\m \, ds =
    \int_{\Omega \times (0,c)} \sqrt{|\nabla f_n^s|^2(x)+|\nabla f_n^x|^2(s)} \, d\m \, ds \\
    & \geq 
    \sup_{\substack{(a,b) \in C(\Omega) \times C(\Omega) \\ a^2+b^2 \leq 1 \\ a,b \geq 0}}
    \int_{\Omega \times (0,c)} a(x) |\nabla f_n^s|(x)+b(x)|\nabla f_n^x|(s) \, d\m \, ds.
    \end{align*}
    Combining with \eqref{E30} and the fact that $f_n=1$ on $\Omega \times \{0\}$ and $f_n=0$ on $\Omega \times \{c\}$, we deduce
    \begin{align*}
        |D f_n|(\Omega \times (0,c)) \geq 
    \sup_{\substack{(a,b) \in C(\Omega) \times C(\Omega) \\ a^2+b^2 \leq 1 \\ a,b \geq 0}}
    \int_\Omega a(x) |\nabla w(f_n)|(x)+b(x) \, d\m=\A(w(f_n),\Omega).
    \end{align*}
    Moreover, $w(f_n) \to w(E)$ in $\sL^1(\Omega)$ as $n \to + \infty$ since
    \[
    \int_\Omega |w(f_n) - w(E)| \, d\m \leq \int_{\Omega \times (0,c)}|f_n(x,s)-1_E(x,s)| \, d\m \, ds.
    \]
    These facts imply that $w(E) \in \BV(\Omega)$ as claimed. Concerning the area of $w(E)$, we get
    \[
    \A(w(E),\Omega) \leq \liminf \A(w(f_n),\Omega) \leq \liminf |D f_n|(\Omega \times (0,c))=\sP(E,\Omega \times \bb{R}),
    \]
    proving the claim \eqref{E|w(E)}.
    \end{proof}
\end{proposition}

If the space in question is finite dimensional, the previous result holds without the assumption that $u$ is bounded. This is the content of Proposition \ref{P|equivalence2}.

\begin{proposition} \label{P|equivalence2}
    Let $(\X,\sd,\m)$ be an $\RCD(K,N)$ space, let $\Omega \subset \X$ be open and let $u \in \BV(\Omega)$, then $u$ is a minimizer of the area functional in $\Omega$ if and only if its hypograph is a perimeter minimizer in $\Omega \times \bb{R}$. 
    \begin{proof}
        If $u$ minimizes the area, then, by Proposition \ref{P|local boundedness}, $u$ is locally bounded in $\Omega$ and we can conclude by Proposition \ref{P|equivalence}. \par 
        Conversely, if the hypograph of $u$ minimizes the perimeter, then $u$ is locally bounded by Proposition \ref{P|aggiunta}. 
        Consider a competitor $f \in \BV(\Omega)$ for $u$ and a subset $\Omega' \subset \subset \Omega$ such that $\{f \neq u\} \subset \subset \Omega'$. We define $f^c:=-c \vee f \wedge c$ for $c:=\|u\|_{\sL^\infty(\Omega')}+1$, so that
        \[
        \A(u,\Omega')=
        \sP(\Hyp(u),\Omega' \times \bb{R})
        \leq \sP(\Hyp(f^c),\Omega' \times \bb{R}) =\A(f^c,\Omega') \leq \A(f,\Omega'),
        \]
        which implies
        $
        \A(u,\Omega) \leq \A(f,\Omega).
        $
    \end{proof}
\end{proposition}

We now prove Theorem \ref{CT4} and its corollary.

\begin{customthm}{4} 
    Let $(\X,\sd,\m)$ be a proper $\RCD(K,\infty)$ space, let $\Omega \subsetneq \X$ be an open set and let $u \in \sL^\infty(\Omega) \cap \W^{1,1}(\Omega)$ be a function whose hypograph is perimeter minimizing.
    If $v \in \sL^\infty(\Omega)$ has hypograph which is perimeter minimizing and $\{u \neq v\} \subset \subset \Omega$, then $u=v$.
    \begin{proof}
    By Theorem \ref{CT1}, $v \in \BV(\Omega)$ and by Proposition \ref{P|equivalence}, both $u$ and $v$ are area minimizers in $\Omega$. 
    
    For every $t \in (0,1)$ we consider the convex combination $tu+(1-t)v$. We have
    \[
    \A(tu+(1-t)v,\Omega) \leq t\A(u,\Omega)+(1-t)\A(v,\Omega)=\A(u,\Omega).
    \]
    Since $tu+(1-t)v$ is a competitor for $u$, the previous inequality implies
    \begin{equation} \label{Ec1}
    \A(tu+(1-t)v,\Omega) =\A(u,\Omega).
    \end{equation}
    The function $a \mapsto \sqrt{1+a^2}$ is strictly convex on the positive real line. Hence, \eqref{Ec1} implies that, for every $t \in (0,1)$, it holds
    \[
    |\nabla(tu+(1-t)v)|=|\nabla u| \quad \m \text{-almost everywhere in } \Omega.
    \]
    Taking the squares in the previous identity and evaluating at $t=0$ and $t=1/2$, we deduce
    \[
    \begin{cases}
        \nabla u \cdot \nabla v = |\nabla u|^2 \\
        |\nabla v|^2=|\nabla u|^2 ,
    \end{cases}
    \]
    which in turn implies 
    \[
    |\nabla(u-v)|=0 \quad \m \text{-almost everywhere in } \Omega.
    \]
    Hence, the function $\psi:=v-u$ satisfies $|D \psi|=|D^s \psi|$ and has support well contained in $\Omega$. Moreover, the function $t \mapsto \A(u+t\psi,\Omega)$ is constant in $(0,1)$. Hence
    \[
    0=\frac{d}{dt} \A(u+t\psi,\Omega) =\frac{d}{dt} |D^s(u+t\psi)|(\Omega)=\frac{d}{dt} |D^s(t\psi)|(\Omega)=|D^s\psi|(\Omega)=|D \psi |(\Omega).
    \]
    In particular, $\psi$ is constant. Since the support of $\psi$ is well contained in $\Omega$ and $\Omega \subsetneq \X$, it holds $\psi \equiv 0$.
    \end{proof}
\end{customthm}

\begin{corollary*}
    Let $(\X,\sd,\m)$ be an $\RCD(K,N)$ space, let $\Omega \subsetneq \X$ be an open set and let $u \in \W^{1,1}(\Omega)$ be a function whose hypograph is perimeter minimizing.
    If $v \in \sL^1(\Omega)$ has hypograph which is perimeter minimizing and $\{u \neq v\} \subset \subset \Omega$, then $u=v$.
    \begin{proof}
        By Lemma \ref{P|aggiunta} both $u$ and $v$ are locally bounded. The statement follows by Theorem \ref{CT4}.
    \end{proof}
\end{corollary*}

\section{Approximation of minimal graphs in \texorpdfstring{$\RCD(K,N)$}{RCD(K,N)}} \label{S|approximation}
The goal of this section is to prove Theorem \ref{T|key approximation}, which implies Theorem \ref{T|approssimazione aggiunta dopo}. \par
We consider spaces $(\X_i,\sd_i,\m_i,x_i)$ converging to a space $(\X,\sd,\m,x)$ in pmGH sense and functions $u_i \in \BV(\X_i)$ converging in $\sL^1$ to $u \in \BV(\X)$. As anticipated, we use the so-called extrinsic viewpoint and we work in the space $(Z,\sd_Z)$ which realizes the convergence. In particular, whenever we consider a ball $B_r(x)$ in this context, unless otherwise specified, we mean the ball of radius $r$ centered in $x$ in the space $Z$ (observe that the intersection of this ball with the embedding of $\X$ in $Z$ would give the corresponding ball in $\X$). 
Similarly, unless otherwise stated, each function $u_i$ will be considered as a function on $Z$ (extending it to zero outside of the embedding of $\X_i$ in $Z$), so that it will make sense
to consider products of the form $u_i 1_{B_r(x)}: Z \to \bb{R}$. To avoid ambiguity, we will denote by $P_t^{\X_i}$ the heat flow on the space $\X_i$.

Similarly, given $f:Z \to \bb{R}$, we will write $\A(f,B_r(x),\m)$ meaning the area of $f$ on $B_r(x) \subset Z$ w.r.t. to the pushforward of the measure $\m$ in $Z$ and we will write $f \in \BV(B_r(x),\m)$ (and the same notation will be used for Lebesgue and Sobolev spaces) to specify that the measure considered is again the pushforward of the measure $\m$ in $Z$. We remark that $f \in \BV(Z,\m)$ if and only if (restricting $f$ to $\X$, where $\m$ is concentrated) $f \in \BV(\X)$, and in this case their total variation measures coincide (see \cite{DiMarinoGigliPratelli}).




\begin{thm} \label{T|local convergence}
    Let $(\X_i,\sd_i,\m_i,x_i)$ be a sequence of $\RCD(K,N)$ spaces converging in pmGH sense to $(\X,\sd,\m,x)$. 
    \begin{enumerate}
    \item If $u \in \BV(\X,\m)$, and $u_i \in \BV(\X_i,\m_i)$ is a sequence such that $u_i \to u$ in $\sL^1$, 
    then for every open set $A \subset Z$ we have
    \[
    \A(u,A,\m) \leq \liminf_{i} \A(u_i,A,\m_i).
    \]
    \item  
    If $u \in \BV(B_r(x),\m)$,
    then for every $s \in (0,r)$ such that $|D^su|(\partial B_s(x))=0$
    there exists a sequence of Lipschitz functions $u_i \in \sL^1(B_s(x),\m_i)$ such that
    $u_i 1_{B_s(x)}$ converges in $\sL^1$-sense to $u 1_{B_s(x)} \in \sL^1(\X,\m)$ and $\A(u_i,B_s(x),\m_i) \to \A(u,B_s(x),\m)$.
    \end{enumerate}
    \begin{proof}
    The first assertion follows by combining the localized lower semicontinuity of total variations given in
    Proposition \ref{P|lsc variation aperti} and the identification between area and perimeter given in Theorem \ref{CT1}. \par 
    We now prove the second assertion following the proof in \cite[Theorem $8.1$]{AHconvergence}. \par 
    Fix $s \in (0,r)$ and observe that, modulo multiplying with a cut off function, we can assume that $u$ has support well contained in $B_r(x)$. Also, modulo truncating (and using a diagonal argument), we can also assume that $u \in \sL^\infty(B_r(x),\m)$ (as the areas of the truncated functions converge to the area of the initial function). 
    In what follows, when we refer to $u$, we implicitly mean its extension to zero on the whole $Z$. \par 
    Since $|D^su|(\partial B_s(x))=0$, Theorem \ref{CT0} guarantees that 
    \[
    \A(P^\X_t(u),B_s(x),\m)
    \to \A(u,B_s(x),\m) \quad
    \text{as } t \to 0.
    \]
    Hence, if we manage to approximate in area the functions $P^\X_t(u)$, we can then conclude with a diagonal argument.
    Each function $P^\X_t(u)$ belongs to $\W^{1,2}(Z,\m)$, so that by Remark \ref{R|asymptotic} we can consider a sequence $f_k \in \Lip_c(Z)$ converging in $\sL^2(Z,\m)$ to $P^\X_t(u)$ and such that $\lip_a(f_k) \to |\nabla P^\X_t(u)|$ in $\sL^2(Z,\m)$.
    We then fix $k$, and as $i$ varies we get 
    \begin{align*}
    \A(f_k,B_{s}(x),\m_i) & \leq
    \int_{B_{s}(x)} \sqrt{1+\lip^2_af_k}\, d\m_i
    \\
    & \leq \int_{Z} 1_{\bar{B}_{s}(x)} \sqrt{1+\lip^2_af_k}\, d\m_i.
    \end{align*}
    Taking a limit as $i$ goes to infinity and taking into account that the last integrand is upper semicontinuous and that $\m(\partial B_{s}(x))=0$, we then obtain
    \[
    \limsup_{i \to + \infty}\A(f_k,B_{s}(x),\m_i) \leq \int_{Z} 1_{\bar{B}_{s}(x)} \sqrt{1+\lip^2_af_k}\, d\m= \int_{B_{s}(x)} \sqrt{1+\lip^2_af_k} \, d\m.
    \]
    Reindexing the sequence $\{f_k\}_{k \in \bb{N}}$ (so that each function is repeated a sufficient number of times) we then get
    \[
    \limsup_{k \to + \infty}\A(f_k,B_{s}(x),\m_{k}) \leq \int_{B_{s}(x)} \sqrt{1+|\nabla P_t(u)|^2} \, d\m = \A(P_t(u),B_{s}(x),\m),
    \]
    which is the desired approximation.
    \end{proof}
\end{thm}

The same proof of the previous theorem can be repeated in the global setting to obtain the $\Gamma$-convergence of the area functional.
In particular, replacing Proposition \ref{P|lsc variation aperti} with the results of \cite{AHconvergence}, we deduce the aforementioned convergence in the case $N=\infty$.

\begin{thm}
    Let $(\X_i,\sd_i,\m_i,x_i)$ be a sequence of proper $\RCD(K,\infty)$ spaces converging in pmGH sense to $(\X,\sd,\m,x)$. 
    \begin{enumerate}
    \item If $u \in \BV(\X,\m)$, and $u_i \in \BV(\X_i,\m_i)$ is a sequence such that $u_i \to u$ in $\sL^1$-sense, then 
    \[
    \A(u,\X,\m) \leq \liminf_{i} \A(u_i,\X_i,\m_i).
    \]
    \item  
    If $u \in \BV(\X,\m)$,
    then
    there exists a sequence of functions $u_i \in \BV(\X_i,\m_i)$ converging in area to $u$.
    \end{enumerate}
    
\end{thm}

We are now in position to prove the main result of the section.

\begin{thm} \label{T|key approximation}
     Let $(\X_i,\sd_i,\m_i,x_i)$ be a sequence of $\RCD(K,N)$ spaces converging in pmGH sense to $(\X,\sd,\m,x)$. Let $r>0$ be such that $B_r(x) \subsetneq \X$ and let $u \in \BV(B_r(x),\m)$ be an area minimizer. 
    For every $s \in (0,r)$, there exists $0<\delta<r-s$ such that the following happens.
    
     There exist area minimizers $\tilde{u}_i \in \BV(B_s(x),\m_i)$ and $\tilde{u} \in \BV(B_{s+\delta}(x),\m)$ such that $\tilde{u}_i 1_{B_s(x)}$ converges in $\sL_1$ to $\tilde{u}1_{B_s(x)}$ and $\{\tilde{u} \neq u\} \subset \subset B_{s+\delta}(x)$. 
     \par 
    Moreover, if $|D^s \tilde{u}|(\partial B_s(x))=0$, it holds $\A(\tilde{u}_i,B_s(x),\m_i) \to \A(\tilde{u},B_s(x),\m)$.
    \begin{proof}
       Let $0<\epsilon<(r-s)/2$ be such that $|D^su|(\partial B_{s+2\epsilon}(x))=0$. We set $\delta:=2\epsilon$. Consider the sequence $u_i \in \BV(B_{s+2\epsilon}(x),\m_i) \cap \sL^{\infty}(B_{s+2\epsilon}(x),\m_i)$ converging in area to $u$ given by Theorem \ref{T|local convergence}. 
       Modulo truncating, we can also assume that $\|u_i\|_{\sL^\infty(B_{s+2\epsilon}(x),\m_i)} \leq \|u\|_{\sL^\infty(B_{s+2\epsilon}(x),\m)}$.
        \par 
        Then, we consider for every $i \in \bb{N}$ the function $\tilde{u}_i \in \BV(B_{s+2\epsilon}(x),\m_i)$ such that $u_i=\tilde{u}_i$ on $B_{s+2\epsilon}(x) \setminus B_{s+\epsilon}(x)$
        and
        \[
        \A(\tilde{u}_i,B_{s+2\epsilon}(x),\m_i) \leq \A(f,B_{s+2\epsilon}(x),\m_i)
        \]
        for every $f \in \BV(B_{s+2\epsilon}(x),\m_i)$ such that $\{f \neq \tilde{u}_i\} \subset B_{s+\epsilon}(x)$.
        Existence of the functions $\tilde{u}_i$ follows by the direct method of calculus of variations. We also note that the condition $u_i=\tilde{u}_i$ on $B_{s+2\epsilon}(x) \setminus B_{s+\epsilon}(x)$ is nonempty by the assumption that $B_r(x) \neq \X$.
        Observe moreover that each $\tilde{u}_i$ satisfies $\|\tilde{u}_i \|_{\sL^\infty(B_{s+2\epsilon}(x),\m_i)} \leq \|u_i \|_{\sL^\infty(B_{s+2\epsilon}(x),\m_i)}$ since, if this is not the case, we can truncate and obtain better competitors. \par 
    Consider then the function $\tilde{u} \in \BV(B_{s+2\epsilon}(x),\m)$ which is the $\sL^1$-limit in $B_{s+2\epsilon}(x)$ of a (non relabeled) subsequence of the functions $\tilde{u}_i$ (this exists thanks to Proposition \ref{P|compactness} and a cutoff argument). We now claim that $\tilde{u}$ is an area minimizer in $B_{s+2\epsilon}(x)$ w.r.t. $\m$. \par 
    To this aim, it is sufficient to show that
        \begin{equation} \label{E|f1}
    \A(\tilde{u},B_{s+2\epsilon}(x),\m) = \A(u,B_{s+2\epsilon}(x),\m).
    \end{equation}
    Observe that by construction $\{u \neq \tilde{u}\} \subset \subset B_{s+2\epsilon}(x)$ and for every $i$ we have $\A(\tilde{u}_i,B_{s+2\epsilon}(x)) \leq \A(u_i,B_{s+2\epsilon}(x))$.
    Hence, since $u$ is an area minimizer, we have that
    \begin{align} \label{E25}
    \nonumber
    \A(\tilde{u},B_{s+2\epsilon}(x),\m) & \geq \A(u,B_{s+2\epsilon}(x),\m)=\lim_{i \to + \infty} \A(u_i,B_{s+2\epsilon}(x),\m_i) \\
    & \geq 
    \limsup_{i \to + \infty} \A(\tilde{u}_i,B_{s+2\epsilon}(x),\m_i)
    \geq \A(\tilde{u},B_{s+2\epsilon}(x),\m),
    \end{align}
    so that
    \[
    \A(\tilde{u},B_{s+2\epsilon}(x),\m) = \A(u,B_{s+2\epsilon}(x),\m).
    \]
    This concludes the proof of the statement. \par
    Suppose now that $|D^s \tilde{u}|(\partial B_s(x))=0$. By \eqref{E25} we know that
    \[
    \A(\tilde{u}_i
    ,B_{s+2\epsilon}
    (x),\m_i) \to \A(\tilde{u},B_{s+2\epsilon}(x),\m).
    \]
    Convergence in area then follows by a standard lower semicontinuity argument combining the fact that $\A(\tilde{u},\partial B_s(x),\m)=0$ and
     point $1$ of Theorem \ref{T|local convergence}.
    \end{proof}
\end{thm}

In the next section, when using the previous theorem, we will occasionally say that the functions $u_i$ converge in $\sL^1$-sense to $u$ in $B_s(x)$, meaning that $u_i1_{B_s(x)} \to u 1_{B_s(x)}$ in $\sL^1$-sense.
Theorem \ref{T|key approximation} immediately implies Theorem \ref{T|approssimazione aggiunta dopo} from the Introduction.

\begin{customthm}{5} 
     Let $(\X_i,\sd_i,\m_i,x_i)$ be a sequence of $\RCD(K,N)$ spaces converging in pmGH sense to $(\X,\sd,\m,x)$. Let $r>0$ be such that $B_r(x) \subsetneq \X$ and let $u \in \W^{1,1}(B_r(x),\m)$ be an area minimizer. 
    For every $s \in (0,r)$, modulo passing to a subsequence, there exist area minimizers $u_i \in \BV(B_s(x),\m_i)$ such that $u_i 1_{B_s(x)}$ converge in $\sL_1$ to $u 1_{B_s(x)}$ and $\A(u_i,B_s(x),\m_i) \to \A(u,B_s(x),\m)$.
    \begin{proof}
        By Theorem \ref{T|key approximation}, there exist area minimizers $u_i \in \BV(B_s(x),\m_i)$ and $\tilde{u} \in \BV(B_{r}(x),\m)$ such that $u_i 1_{B_s(x)}$ converges in $\sL_1$ to $\tilde{u}1_{B_s(x)}$ and $\{\tilde{u} \neq u\} \subset \subset B_{r}(x)$. By the corollary of Theorem \ref{CT4}, it holds $\tilde{u}=u$. By the last part of Theorem \ref{T|key approximation}, since $\m(\partial B_s(x))=0$ for every $s>0$, it holds $\A(u_i,B_s(x),\m_i) \to \A(u,B_s(x),\m)$.
    \end{proof}
    \end{customthm}

\appendix
\section{Appendix}



We collect here some applications of the previously obtained theorems in the setting of Ricci limit spaces. 
Here, we work under the additional assumption that Riemannian manifolds with a lower bound on the Ricci curvature satisfy a mean value property for minimal graphs, i.e. assumption \eqref{Assumption} below.
We recall some notation. Given a manifold $(\M,g)$, a solution of the minimal surface equation $u \in C^\infty(\M)$, and a point $p \in \M \times \bb{R}$, we denote by $B^\times_r(p) \subset \M \times \bb{R}$ a ball of radius $r$ centered in $p$ w.r.t. the product distance in $\M \times \bb{R}$, and we set
\[
B^m_r(p):=B^\times_r(p) \cap \graph(u) \subset \M \times \bb{R}.
\]
Given a Riemannin manifold $(\M,g)$, we denote by $\aH^n$ the $n$-dimensional Hausdorff measure induced by the Riemannian distance.
\begin{assumption}{MVP} \label{Assumption}
Let $(\M^n,g)$ be a manifold with $\Ric_\M \geq K$ and $\aH^n(B_1(x)) \geq v>0$. There exists $c(K,n,v)>1$ satisfying the following. Let $u \in C^\infty(B_{c}(x))$ be a solution of the minimal surface equation. 
    If $f_1,f_2 \in C^\infty(B^m_c(x,u(x)))$ are non-negative functions which are respectively sub-harmonic and super-harmonic, then
    \[
    \sup_{B^m_{1/2}(x,u(x))} f_1 \leq c\fint_{B^m_{1}(x,u(x))} f_1 \, d \aH^n \quad \text{and} \quad
    \fint_{B^m_{1}(x,u(x))} f_2 \, d \aH^n
    \leq c \inf_{B^m_{1/2}(x,u(x))} f_2.
    \]
\end{assumption}

We spend a few words on the assumption \eqref{Assumption}.
The mean value property \eqref{Assumption} would follow if one is able to prove Sobolev and Poincaré-Neumann inequalities for functions defined on minimal graphs on Riemannian manifolds with Ricci curvature bounded from below. This is the strategy used in \cite{DP23} to show that \eqref{Assumption} holds. Unfortunately, some parts of the proofs of the aforementioned inequalities in \cite{DP23} are unclear. For this reason, we add the mean value property as an assumption. 
Nevertheless, it is reasonable to think that these functional inequalities on manifolds should hold. 

The main results of this part of the paper are Theorems \ref{A1} and \ref{A2} below.
Theorem \ref{A1} shows that if \eqref{Assumption} holds, Sobolev minimizers of the area functional on non-collapsed Ricci limit spaces are locally Lipschitz. This is remarkable as minimizers of the area on generic metric measure spaces are much less regular, as shown in \cite{HKLreg} (see also Example \ref{E|new} for the existence of a Sobolev minimizer of the area on a 'bad' metric space which is not Lipschitz).

\begin{thm} \label{A1}
Assume that \eqref{Assumption} holds. Let $v>0$, $n \in \bb{N}$ and $K \in \bb{R}$ be fixed.
    Let $(\M^n_k,g_k,x_k)$ be a sequence of pointed Riemannian manifolds with $\Ric_{\M^n_k} \geq K$ and $\aH^n_k(B_1(x_k)) \geq v >0$ converging in pointed Gromov-Hausdorff sense to $(\X,\sd,x)$. If $\Omega \subset \X$ is open and $u \in \W^{1,1}(\Omega,\aH^n)$ is an area minimizer, then $u \in \Lip_{loc}(\Omega)$.
\end{thm}

Theorem \ref{A1} is proved combining the a-priori gradient estimates for minimizers of the area on manifolds obtained in \cite{DP23} (which are a consequence of \eqref{Assumption}) and Theorem \ref{T|approssimazione aggiunta dopo}. These gradient estimates are also shown to hold for Ricci limit spaces (see Theorem \ref{T|previously CT6} and the stronger version given by Theorem \ref{T|previously CT9} in the presence of maximal volume growth).

Finally, Theorem \ref{A2} is part of a series of recent results proving Bernstein-type theorems on manifolds with a lower Ricci curvature bound (see \cite{DingJost,Ding,CMMR,cmmr23,DingCap,DP23}) and on metric measure spaces with with an analogous synthetic bound (see \cite{C1,C2}).
The proof, once again, relies on the one of the analogous result on Riemannian manifolds proved in \cite{DP23} and Theorem \ref{T|approssimazione aggiunta dopo}.

\begin{thm} \label{A2}
    Assume that \eqref{Assumption} holds. Let $n \in \bb{N}$ be fixed and let $(\M^n_k,g_k,x_k)$ be a sequence of pointed Riemannian manifolds with $\Ric_{\M^n_k} \geq 0$ converging in pointed Gromov-Hausdorff sense to a space $(\X,\sd,x)$ such that
    \[
    \lim_{r \to + \infty} \frac{\aH^n(B_r(x))}{r^n} \geq v >0.
    \]
    If $u \in \W^{1,1}_{loc}(\X,\aH^n)$ is an area minimizer such that
    \[
    \lim_{y \to +\infty} \frac{\max \{-u(y),0\}}{\sd(x,y)}=0,
    \]
    then $u$ is constant.
\end{thm}

Before proving Theorems \ref{A1} and \ref{A2}, we recall some notation.
Given a metric measure space $(\X,\sd,\m)$, a function $u \in \Lip_{loc}(\X)$, and $(x,t) \in \X \times \bb{R}$, we will use the following notation for balls on the graph of $u$ (where the '$m$' stands for minimal, as $u$ will be a solution of the minimal surface equation):
\[
B^m_r(x,t):=B^\times_r(x,t) \cap \graph(u) \subset \X \times \bb{R}.
\]
Similarly, given a point $x \in \X$, we will denote by $\bar{x}$ its projection on the graph of $u$, i.e. $\bar{x}:=(x,u(x))$. Given a metric space $(\X,\sd)$, we denote by $\aH^n$ the $n$-dimensional Hausdorff measure induced by $\sd$.
\par 
In this section, we will often work in the setting of Riemannian manifolds. Whenever we say that we consider a Riemannian manifold $(\M^n,\sd)$, we implicitly mean that we have a Riemannian metric $g$ on $\M$ and that $\sd$ is the Riemannian distance. 
Given a manifold $(\M,g)$ and a solution of the minimal surface equation $u \in C^\infty(\Omega)$ on an open subset $\Omega \subset \M$, we denote by $\nabla_m$ and $\Delta_m$ respectively the gradient and the Laplacian on the graph of $u$. 
In the next proposition and more in general in the paper we will consider the laplacian $\Delta_m$ of a function $f$ defined on $\Omega \subset \M$. When doing so, we implicitly mean that we are applying the laplacian to the function $f \circ \pi$, where $\pi:\graph(u) \to \Omega$ is the standard projection. We will often use the fact that minimizers of the area functional on smooth manifolds are smooth. 

We now turn our attention to proving Theorems \ref{A1} and \ref{A2}.
The next proposition is the so-called Jacobi equation for minimal graphs and it follows by taking the second variation of the area of a minimal graph.

\begin{proposition} \label{P|superharmonicity}
     Let $(\M,g)$ be a Riemannian manifold, let $\Omega \subset \M$ be open and let $u \in C^\infty( \Omega)$ be a solution of the minimal surface equation. Then, setting $v_u:=\sqrt{1+|\nabla u|^2}$ and denoting by $A$ the second fundamental form of $\graph(u)$ in $\Omega \times \bb{R}$, we have
     \[
     \Delta_m v_u^{-1} = -\frac{|A|^2}{v_u}-\frac{\Ric(\nabla u,\nabla u)}{v_u^{3/2}},
     \]
     or, equivalently,
     \[
     \Delta_m \log v_u= |A|^2 + \frac{\Ric(\nabla u,\nabla u)}{v_u^2}+|\nabla_m \log v_u|^2.
     \]
\end{proposition}

The mean value property \eqref{Assumption}, combined with Proposition \ref{P|superharmonicity}, gives a-priori gradient estimates for minimal graphs. 
This is the content of Theorem \ref{T|Ding} (see \cite[Theorem $6.1$]{DP23}).
Although stated in a slightly different form, inspecting the proof of \cite[Theorem $6.1$]{DP23} (and using \eqref{Assumption} instead of \cite[Lemma 5.7]{DP23}), one can check that this is the version we obtain.


\begin{thm} \label{T|Ding}
    Assume that \eqref{Assumption} holds. Let $r>0$ and let $(\M^n,g)$ be a complete non-compact Riemannian manifold
    with 
    \[
    \Ric_{\M} \geq -(n-1)k^2r^{-2}
    \]
    on $B_r(p) \subset \M$ for some $k \geq 0$. Suppose that $\aH^n(B_r(p)) \geq v r^n$ for some constant $v>0$. There exists $c(k,n,v)>1$ such that if $u \in C^\infty(B_{cr}(p))$ is a solution of the minimal surface equation, then
    \[
    \sup_{{x \in B^m_{r}(\bar{p})}}
    |\nabla u(x)| \leq c e^{cr^{-1}(u(p)-\inf_{x \in B_{cr}(p)}u(x))}.
    \]
\end{thm}

We now prove three basic lemmas.

\begin{lemma} \label{L|noncollapsing limit}
    Let $n \in \bb{N}$ and $K \in \bb{R}$ be fixed. Let $(\M^n_k,\sd_k,x_k)$ be a sequence of pointed Riemannian manifolds with $\Ric_{\M^n_k} \geq K$ converging in pGH sense to $(\X,\sd,x)$. For every $r>0$, if $k$ is large enough, then $\aH^n_k(B_r(x_k)) \geq \aH^n(B_r(x))/2$.
    \begin{proof}
        If $\aH^n(B_r(x))=0$ the statement follows trivially. If instead $\aH^n(B_r(x))>0$, then the Hausdorff dimension of $(\X,\sd)$ is at least $n$, so that Theorem \ref{T|Gigli De philippis} gives that the sequence $(\M^n_k,\sd_k,\aH^n_k,x_k)$ converges in pmGH sense to $(\X^n,\sd,\aH^n,x)$. The statement then follows by lower semicontinuity of measures w.r.t. weak convergence.
    \end{proof}
\end{lemma}

In the next lemma, we consider the behavior of the local Lipschitz constant of an area minimizing function $u \in \BV(\Omega)$. Here, we implicitly mean that we are considering a precise representative for $u$, i.e. the function whose value in $x$ is given by
\[
\limsup_{s \to 0} \fint_{B_s(x)} u \, d\m.
\]
This will be implicitly used when we exploit the fact that $\graph(u) \subset \partial \Hyp(u)$, where $\Hyp(u)$ denotes the closed representative of the hypograph of $u$.

\begin{lemma} \label{L|uniform lip}
Let $n \in \bb{N}$ and $K \in \bb{R}$ be fixed. Let $(\M^n_i,\sd_i,\aH^n_i,x_i)$ be a sequence of pointed Riemannian manifolds with $\Ric_{\M^n_i} \geq K$ converging in pmGH sense to $(\X,\sd,\m,x)$. If $u \in \BV(B_r(x),\m)$ is an area minimizer and there is a sequence
of area minimizers $u_i \in \BV(B_R(x_i),\aH^n_i)$ with $R>r$ converging in $\sL_1$ to $u$ on $B_r(x)$ with
\[
\sup_{B_R(x_i)} \lip (u_i) \leq c,
\]
then
\[
\sup_{B_r(x)} \lip( u) \leq c.
\]
\begin{proof}
    Suppose by contradiction that the statement is false. Therefore, there exists $y \in B_r(x) \cap \X$ and a sequence $y_k \to y$ contained in $B_r(x) \cap \X$ such that
\[
\frac{|u(y)-u(y_k)|}{\sd(y,y_k)} > c .
\]
Since the boundaries of the hypographs of $u_i$ converge in Kuratowski sense in $B_r(x) \times \bb{R}$ (as a subset of the space $(Z,\sd_z)$ realizing the pmGH convergence) to $\partial \Hyp(u)$, and since $\graph(u) \subset \partial \Hyp(u)$, we have that for every $k$ there exists a sequence $y_i^k \in B_r(x) \cap \M_i$ such that $y_i^k \to y_k$ in $Z$ and $u_i(y_i^k) \to u(y_k)$.
Similarly, there exists $y_i \in B_r(x) \cap \M_i$ such that $y_i \to y$ in $Z$ and $u_i(y_i) \to u(y)$. This implies that for every $k>0$ fixed, for $i$ large enough, we have
\[
\frac{|u_i(y_i)-u_i(y^k_i)|}{\sd_i(y_i,y^k_i)} > c.
\]
In particular, if $k$ and $i$ are sufficiently large, we will have that the geodesic in $\M_i$ connecting $y_i$ and $y^k_i$ is contained in $B_R(x_i)$. Integrating along this geodesic we have
\[
|u_i(y_i)-u_i(y^k_i)| \leq c \sd_i(y_i,y^k_i),
\]
a contradiction.
\end{proof}
\end{lemma}

\begin{lemma} \label{L|projection}
Let $(\M,g)$ be a Riemannian manifold and let $u \in C^{\infty}(B_r(x))$ be a smooth function.
If $0<t \leq r$, $c>1$ and
    \[
    \sup_{B_t^{m}(\bar{x})} \lip(u) \leq c,
    \]
    then
    \[
    \pi(\bar{B}_t^{m}(\bar{x})) \supset B_{t/(2c)}(x).
    \]
    \begin{proof}
        Let $y \in B_{t/(2c)}(x)$ and let $\gamma$ be a length minimizing geodesic connecting $x$ and $y$. If $u$ restricted to this geodesic takes value in $B_{t/(2c)}(x) \times (u(x)-t/2,u(x)+t/2)$, then there is nothing left to prove.
        So suppose by contradiction that this is not true.
        Without loss of generality this implies that that there exists $\gamma(t_0)$ such that $u(\gamma(t_0))=u(x)+t/2$, and we can also suppose that for every $s \in (0,t_0)$ we have $u(\gamma(s)) \in (u(x)-t/2,u(x)+t/2)$. In this case we have that
        \[
        |u(x)-u(\gamma(t_0))| \leq c \sd(\gamma(t_0),x) < t/2,
        \]
        a contradiction.
    \end{proof}
\end{lemma}

The next theorem complements Theorem \ref{A1} and, assuming \eqref{Assumption}, gives explicit a-priori gradient estimates for minimizers of the area on Ricci limit spaces.

\begin{thm} \label{T|previously CT6}
    Assume that \eqref{Assumption} holds. Let $r>0$, $v>0$, $n \in \bb{N}$ and $K \in \bb{R}_+$ be fixed. Let $(\M^n_k,\sd_k,x_k)$ be a sequence of pointed Riemannian manifolds with $\Ric_{\M^n_k} \geq -(n-1)K^2r^{-2}$ and $\aH^n_k(B_r(x_k)) \geq v r^n$ converging in pointed Gromov-Hausdorff sense to $(\X,\sd,x)$ with $B_r(x) \subsetneq \X$. 
    Let $u \in \W^{1,1}(B_r(x),\aH^n)$ be an area minimizer.
    There exists $c(K,n,v)>0$ such that, setting 
    \[
    c_r:=c e^{cr^{-1}(u(x)-\inf_{y \in B_r(x)}u(y))},
    \]
    it holds
     \[
    \sup_{z \in B_{r/c_r}(x)}\lip(u)(z) \leq
     c_r.
    \]
    \begin{proof}
    Fix $r/2<t<s<r$ and let $u_k \in \BV(B_s(x),\aH^n_k)$ be the area minimizers given by Theorem \ref{T|approssimazione aggiunta dopo} such that $u_k1_{B_s(x)} \to u1_{B_s(x)}$ in $\sL^1$ sense.
    Each function $u_k$ is defined on $B_t(x_k)$ for $k$ large enough and it is smooth on this ball.
    Moreover, thanks to Theorem \ref{T|Ding}, for some $c(k,n,v)>1$,we have 
    \begin{equation} \label{E14}
    \sup_{B^m_{t/c}(\bar{x}_k)}|\nabla u_k| \leq c e^{ct^{-1}(u_k(x_k)-\inf_{y \in B_t(x_k)}u_k(y))}.
    \end{equation}
    Observe now that the Kuratowski convergence of the graphs of $u_k$ to the boundary of the hypograph of $u$ in $B_s(x)$ implies that for every $\delta>0$, if $k$ is large enough, we have
    \[
    \inf_{z \in B_{t}(x_k)} u_k(z) \geq \inf_{z \in B_{s}(x)} u(z) - \delta
    \]
    and
    \[
    u_k(x_k) \leq \sup_{y \in B_\delta(x)}u(y)+\delta.
    \]
    Hence, for every $\delta>0$ and for $k$ large enough,
    the inequality \eqref{E14} becomes
    \begin{equation} \label{E13}
    \sup_{B^m_{t/c}(\bar{x}_k)}|\nabla u_k| \leq c e^{ct^{-1}(\sup_{y \in B_\delta(x)}u(y)-\inf_{y \in B_s(x)}u(y)+2\delta)}.
    \end{equation}
    We now set $c^t_\delta:=  c e^{ct^{-1}(\sup_{y \in B_\delta(x)}u(y)-\inf_{y \in B_s(x)}u(y)+2\delta)}$ and we note that, by Lemma \ref{L|projection}, it holds
    \begin{equation*} 
    \sup_{B_{t/(2cc^t_\delta)}(x_k)}|\nabla u_k| \leq c^t_\delta,
    \end{equation*}
    which, together with
    Lemma \ref{L|uniform lip}, implies that
    \begin{equation} \label{E15}
   \sup_{z \in B_{t/(3cc^t_\delta)}(x)} \lip(u) \leq
c^t_\delta .
    \end{equation}
    This shows in particular that $u$ is continuous in $x$, so that when we let $\delta \to 0$ we have that
    \[
    \sup_{y \in B_\delta(x)}u(y) \to u(x).
    \]
    Hence, taking the limit as $\delta \to 0$ in \eqref{E15} and setting $c^t_0:=c e^{ct^{-1}(u(x)-\inf_{y \in B_s(x)}u(y))}$, we get
    \[
    \sup_{z \in B_{t/(3cc^t_0)}(x)} \lip(u)(z) \leq
     c^t_0 \leq c_r:=ce^{2cr^{-1}(u(x)-\inf_{y \in B_r(x)}u(y))}.
    \]
    Now, taking the limit as $t \to r$ and using the same notation, we get
    \[
    \sup_{z \in B_{r/(3cc_r)}(x)} \lip(u)(z) \leq
     c_r,
    \]
    as claimed.
    \end{proof}
\end{thm}

\begin{customthm}{A.1}
Assume that \eqref{Assumption} holds. Let $v>0$, $n \in \bb{N}$ and $K \in \bb{R}$ be fixed.
    Let $(\M^n_k,\sd_k,x_k)$ be a sequence of pointed Riemannian manifolds with $\Ric_{\M^n_k} \geq K$ and $\aH^n_k(B_1(x_k)) \geq v >0$ converging in pointed Gromov-Hausdorff sense to $(\X,\sd,x)$. If $\Omega \subset \X$ is open and $u \in \W^{1,1}(\Omega,\aH^n)$ is an area minimizer, then $u \in \Lip_{loc}(\Omega)$.
    \begin{proof}
        For every $y \in \Omega$ and every sequence $y_k \in \M_k$ converging to $y$, we can find $\epsilon>0$ and $v_y>0$ such that $\aH_k^n(B_\epsilon(y_k)) \geq v_y>0$ for every $k$ large enough (this follows from the doubling property of each of these manifolds and the uniform lower bound on the volume of the unit ball).
        Hence, the previous theorem guarantees that $u \in \Lip(B_\delta(y))$ for $\delta>0$ sufficiently small. By the arbitrariness of $y \in \Omega$, we conclude.
    \end{proof}
\end{customthm}

\begin{example} \label{E|new}
    On generic metric measure spaces, Sobolev minimizers of the area functional are not locally Lipschitz. To see this, consider the metric measure space $(\bb{R},\sd_e,\mu)$, where $\sd_e$ is the Euclidean distance and $\mu:=|x|\lambda^1$. We claim that the function $1_{\bb{R}_+}$ belongs to $\W^{1,1}_{loc}(\bb{R},\sd_e,\mu)$ and minimizes the area.
    Indeed, consider the Lipschitz functions $f_n(x):=(1-n\sd_e(x,\bb{R}_+) \vee 0$. For every $r>0$, it holds
    \[
    |D 1_{\bb{R}_+}|(-r,r) \leq \lim_{n \to + \infty} \int_{-r}^r \lip f_n(x) \, d\mu=
    \lim_{n \to + \infty} n\int_{-1/n}^{0} |x| \, dx=0.
    \]
    This proves that $|D 1_{\bb{R}_+}|=0$, so that, in particular, $1_{\bb{R}_+} \in \W^{1,1}_{loc}(\bb{R},\sd_e,\mu)$ and $\A(1_{\bb{R}_+},\Omega)=\mu(\Omega)$ for every $\Omega \subset \bb{R}$.
    Let now $u \in \BV(\bb{R},\sd_e,\mu)$ be a competitor for $1_{\bb{R}_+}$ with $\{u \neq 1_{\bb{R}_+}\} \subset \subset \Omega$ for some open set $\Omega \subset \subset \bb{R}$. Then, by definition of area, it holds
    \[
    \A(u,\Omega) \geq \mu(\Omega)=\A(1_{\bb{R}_+},\Omega),
    \]
    showing that $1_{\bb{R}_+}$ is indeed a minimizer.
\end{example}
The next lemma, which will be used to prove Theorem \ref{T|previously CT9}, follows by repeating the argument of \cite[Theorem $6.2$]{DP23}, and for this reason some details in the proof are omitted.

\begin{lemma} \label{L|partial global Lip}
    Assume that \eqref{Assumption} holds. Let $r_0>0$, $v>0$, $n \in \bb{N}$. Let $(\M^n,\sd,x)$ be a manifold with $\Ric_{\M^n} \geq 0$ and $\aH^n(B_{r_0}(x)) \geq v r_0^n$. There exists $c(n,v)>1$ such that, if $u \in C^\infty(B_{cr_0}(x))$ is an area minimizer, then given $r>0$ such that
    \[
    r^2 = r_0^2 -4 \Big(\sup_{y \in B_{r/2}(x)}|u(y)-u(x)| \Big)^2,
    \]
    we have
    \[
    \sup_{y \in B_{r/2}(x)} |\nabla u|(y) \leq c\Big(1+r^{-1} \sup_{y \in B_{r/2}(x)}|u(y)-u(x)|\Big)^n.
    \]
    \begin{proof}
        Thanks to Propositions \ref{P|superharmonicity} and \eqref{Assumption},
        there exists a constant $c(n,v)>0$ such that
        \[
        \frac{1}{\aH^{n}( B^m_{r_0}(\bar{x}))} \int_{B^m_{r_0}(\bar{x})} \frac{1}{\sqrt{1+|\nabla u |^2}} \, d\aH^{n} 
        \leq c\inf_{y \in \pi( B^m_{r_0/2}(\bar{x}))} \frac{1}{\sqrt{1+|\nabla u |^2}}(y),
        \]
        where $\aH^n$ is the Hausdorff measure in $\M \times \bb{R}$ and we are using the previously defined convention for functions defined on $\M$ evaluated on the graph of $u$. \par 
        Hence, possibly changing $c$ at every step, we get
        \[
        \sup_{y \in \pi( B^m_{r_0/2}(\bar{x}))} \sqrt{1+|\nabla u |^2}(y) \leq c \frac{\aH^{n}( B^m_{r_0}(\bar{x}))}{\aH^{n}(\pi(  B^m_{r_0}(\bar{x})))},
        \]
        which, together with the fact that $u$ is an area minimizer, implies that
        \[
        \sup_{y \in \pi(  B^m_{r_0/2}(\bar{x}))} \sqrt{1+|\nabla u |^2}(y) \leq \frac{c r_0^n}{\aH^{n}(\pi(  B^m_{r_0}(\bar{x})))}.
        \]
        By our hypothesis on $r>0$ we then get that $B_{r/2}(x) \subset \pi(  B^m_{r_0/2}(\bar{x})) \subset \pi(  B^m_{r_0}(\bar{x}))$, while Bishop Gromov's inequality and the fact that 
        $\aH^n(B_{r_0}(x)) \geq v r_0^n$ imply that $\aH^n(B_{r/2}(x)) \geq 2^{-n}v r^n$. Putting these facts together we get, changing again $c>0$,
        \[
         \sup_{y \in \pi(B^m_{r_0/2}(\bar{x}))} \sqrt{1+|\nabla u |^2}(y) \leq c \frac{ r_0^n}{r^n}
         \leq c \Big(1+r^{-2} \sup_{y \in B_{r/2}(x)}|u(y)-u(x)|^2\Big)^{n/2},
        \]
        which then implies the statement.
    \end{proof}
\end{lemma}

The next theorem, assuming \eqref{Assumption}, gives improved a-priori gradient estimates for minimizers of the area in Ricci limit spaces with maximal volume growth. The proof is adapted from \cite[Theorem 6.2]{DP23}.

\begin{thm} \label{T|previously CT9}
    Assume that \eqref{Assumption} holds. Let $n \in \bb{N}$ be fixed and let $(\M^n_k,\sd_k,x_k)$ be a sequence of pointed Riemannian manifolds with $\Ric_{\M^n_k} \geq 0$ converging in pointed Gromov-Hausdorff sense to a space $(\X,\sd,x)$ such that
    \[
    \lim_{r \to + \infty} \frac{\aH^n(B_r(x))}{r^n} \geq v >0.
    \]
    If $u \in \W^{1,1}_{loc}(\X,\aH^n)$ is an area minimizer, then there exists $c(n,v)>0$ such that, for every $r>0$, it holds
    \[
    \sup_{y \in B_{r/2}(x)} \lip( u)(y) \leq c\Big(1+r^{-1} \sup_{y \in B_r(x)}|u(y)-u(x)|\Big)^n.
    \]    
    \begin{proof}
        Fix $0<R_1<R_2$ and consider the sequence of area minimizers $u_k \in \BV(B_{R_2}(x),\aH^n_k)$ converging in $\sL^1$ to $u$ on $B_{R_2}(x)$ given by Theorem \ref{T|approssimazione aggiunta dopo}. 
    Each function $u_k$ is defined on $B_{R_1}(x_k)$ for $k$ large enough and is smooth on this ball. \par
    Consider then $r>0$ such that
    \[
    r^2 = R_1^2 -4 \Big(\sup_{y \in B_{r/2}(x)}|u(y)-u(x)| \Big)^2.
    \]
    Since $u$ is locally Lipschitz by Theorem \ref{A2}, the Kuratowski convergence of the hypographs of $u_k$ to the hypograph of $u$ implies that
    \[ 
    \sup_{y \in B_{r/2}(x_k)}|u_k(y)-u_k(x_k)| \to \sup_{y \in B_{r/2}(x)}|u(y)-u(x)|.
    \]
    In particular, setting $r_k >0$ to be such that
    \[
    r_k^2 = R_1^2 -4\Big(\sup_{y \in B_{r_k/2}(x_k)}|u_k(y)-u(x_k)| \Big)^2,
    \]
    we have that $r_k \to r$ as $k \to \infty$. \par 
    Moreover, combining Bishop Gromov's inequality on $\X$ with Lemma \ref{L|noncollapsing limit}, we obtain that, for every $k$ large enough, we have
    \[
    \aH^n_k(B_{R_1}(x_k)) \geq \frac{v}{2}R_1^n,
    \]
    so that, by Lemma \ref{L|partial global Lip}, we obtain
    \[
    \sup_{y \in B_{r_k/2}(x)} |\nabla u_k|(y) \leq c(n,v)\Big(1+{r_k}^{-1} \sup_{y \in B_{r_k/2}(x_k)}|u_k(y)-u_k(x_k)|\Big)^n.
    \]
    Moreover, using Lemma \ref{L|uniform lip} and the fact that $r_i \to r$, we get 
    \begin{equation} \label{E16}
    \sup_{y \in B_{r/4}(x)} |\nabla u|(y) \leq c\Big(1+{r}^{-1} \sup_{y \in B_{r/2}(x)}|u(y)-u(x)|\Big)^n.
    \end{equation}
    Finally, the function $R_1 \mapsto r$ is increasing (by definition), continuous (since $u$ is continuous), and tends to $+\infty$ at infinity (since $u$ is locally bounded), so that \eqref{E16} holds for every $r>0$ by the arbitrariness of $R_1$.
    \end{proof}
\end{thm}

The next lemma will be used to prove Theorem \ref{A2}. The proof is adapted from the one in \cite[beginning of Theorem $3.6$]{DingJost}.

\begin{lemma}
    Let $(\X,\sd,\m,x)$ be an $\RCD(0,N)$ space.
    If $u \in \Lip(\X)$ is a function whose hypograph minimizes the perimeter and such that
    \begin{equation} \label{E20}
    \lim_{y \to +\infty} \frac{\max \{-u(y),0\}}{\sd(x,y)}=0,
    \end{equation}
    then
    \[
    \lim_{y \to +\infty} \frac{|u(y)|}{\sd(x,y)}=0.
    \]
    \begin{proof}
        Let $L:=\ssf{L}(u) \vee 1$, where $\ssf{L}(u)$ is the Lipschitz constant of $u$. The Harnack inequality on the graph of $u$ (Theorem \ref{T|Harnack}) gives that, for every $r>0$, for some $C>1$, it holds
        \[
        \sup_{y \in B^m_r(\bar{x})} u(y) - \inf_{z \in B^m_r(\bar{x})} u(z) \leq C(u(x)-\inf_{z \in B^m_r(\bar{x})} u(z)),
        \]
        so that
        \[
        \sup_{y \in B^m_r(\bar{x})} u(y) \leq Cu(x)-(C-1)\inf_{z \in B^m_r(\bar{x})} u(z).
        \]
        Using that $u$ is Lipschitz, we get
        \[
        \sup_{y \in B_{r/(2L)}(x)} u(y) \leq Cu(x)-(C-1)\inf_{z \in B_{r}(x)} u(z).
        \]
        Moreover,
        from \eqref{E20}, for every $\delta>0$ there exists $C_\delta>0$ such that $u(y) \geq -C_\delta-\delta \sd(x,y)$ for every $y \in \X$. Hence,
        \[
        \sup_{y \in B_{r/(2L)}(x)} u(y) \leq Cu(x)+(C-1)(C_\delta+\delta r).
        \]
        This implies that
        \[
    \limsup_{y \to +\infty} \frac{|u(y)|}{\sd(x,y)} \leq 2CL\delta,
    \]
    which implies the statement by the arbitrariness of $\delta>0$.
    \end{proof}
\end{lemma}

The proof of the next result, i.e. Theorem \ref{A2}, is adapted from the one in \cite[Theorem $3.6$]{DingJost}.

\begin{customthm}{A.2}
    Assume that \eqref{Assumption} holds. Let $n \in \bb{N}$ be fixed and let $(\M^n_k,\sd_k,x_k)$ be a sequence of pointed Riemannian manifolds with $\Ric_{\M^n_k} \geq 0$ converging in pointed Gromov-Hausdorff sense to a space $(\X,\sd,x)$ such that
    \[
    \lim_{r \to + \infty} \frac{\aH^n(B_r(x))}{r^n} \geq v >0.
    \]
    If $u \in \W^{1,1}_{loc}(\X,\aH^n)$ is an area minimizer such that
    \[
    \lim_{y \to +\infty} \frac{\max \{-u(y),0\}}{\sd(x,y)}=0,
    \]
    then $u$ is constant.
    \begin{proof}
        By Lemma \ref{L|noncollapsing limit} and Theorem \ref{T|previously CT6} we know that $u \in \Lip(\X)$. Hence the previous lemma gives that
        \[
        \lim_{y \to +\infty} \frac{|u(y)|}{\sd(x,y)}=0.
        \]
        We claim that $|\nabla u|=0$ in $B_1(x)$. A rescaling argument then gives $|\nabla u|=0$ in $B_r(x)$ for every $r>0$, yielding the statement. \par 
        To this aim, we fix $\epsilon>0$ and we choose $R>1$ large enough so that, for every $y \in \X$ with $\sd(x,y)>R$, it holds
        \[
        |u(y)| \leq \epsilon \sd(x,y).
        \]
        Consider $R_1>R$ to be fixed later and the minimizers $u_k \in C^\infty(B_{5R_1}(x) \cap \M_k)$ converging in $\sL^1$ to $u$ on $B_{5R_1}(x)$.
        For every $k$ large enough, since the hypographs of $u_k$ converge in Kuratowski sense to the hypograph of $u$ in $B_{5R_1}(x) \times \bb{R}$, we have that
         \begin{equation} \label{E23}
        |u_k(y)| \leq 2\epsilon \sd_k(x,y)
        \end{equation}
        for every $y \in B_{R_1}(x_k)$ such that $\sd_k(x_k,y)>2R$. \par
        So far we estimated the growth of the functions $u_k$ exploiting the growth condition on $u$; now we estimate their Lipschitz constants. 
        For every $s>0$ let $c_s$ be the quantity defined in Theorem \ref{T|previously CT6} relative to $u$ and let $c_s^k$ denote the same quantity relative to each $u_k$. 
        Note that by our assumption on the growth of $u$, we have that there exists $c>1$ such that $c_s \leq c$ for every $s$.\par 
        Since the hypographs of the functions $u_k$ converge in Kuratowski sense to the one of $u$, we have that $c_{4R_1}^k \to c_{4R_1}$, so that by the a-priori gradient estimate Theorem \ref{T|previously CT6}, for $k$ large enough,
        \begin{equation} \label{E22}
    \sup_{z \in B_{3R_1/(16c)}(x_k)}|\nabla u_k|(z) \leq
     2c.
    \end{equation}
    To shorten the notation we set $R_2:=3R_1/(16c)$.
    A posteriori, since \eqref{E22} holds for every $R_1$ large enough, we see that we could have chosen $R_1>0$ (large enough) so that
        \begin{equation} \label{E21}
        B_{2R}(x_k) \subset \pi(B^m_{R_2}(\bar{x}_k)).
        \end{equation}
        We now denote by $\aH^n_k$ and $\nabla_\times$ respectively the Hausdorff measure and the gradient in $\M_k \times \bb{R}$ w.r.t. the product distance, and by $\nabla_m$ the gradient on the graph of $u_k$. We claim that
        \begin{equation} \label{E24}
        \int_{B^m_{R_2}(\bar{x}_k)}|\nabla_m u_k|^2 \, d \aH^n_k \leq 4 \epsilon^2 \aH^n_k(B^m_{2R_2}(\bar{x}_k)).
        \end{equation}
        To prove the claim, we consider $\eta \in \Lip(B^\times_{2R_2}(\bar{x}_k))$ with $\eta=0$ on $\partial B^\times_{2R_2}(\bar{x}_k)$, $\eta=1$ on $B^\times_{R_2}(\bar{x}_k)$ and $\nabla_\times \eta \leq 1/R_2$. Using the fact that $u_k$ is harmonic on its graph, we then get
        \[
        \int_{B^m_{R_2}(\bar{x}_k)}|\nabla_m u_k|^2 \, d \aH^n_k
        \leq 
        \int_{B^m_{2R_2}(\bar{x}_k)}|\nabla_m u_k|^2 \eta \, d \aH^n_k
        \leq 4\int_{B^m_{2R_2}(\bar{x}_k)} u_k^2 |\nabla_m \eta|^2 \, d \aH^n_k. 
        \]
        Moreover, combining \eqref{E23}, \eqref{E21} and the fact that $|\nabla_m \eta| \leq |\nabla_\times \eta|$, we deduce that
        \[
         4\int_{B^m_{2R_2}(\bar{x}_k)} u_k^2 |\nabla_m \eta|^2 \, d \aH^n_k 
        \leq 4 \epsilon^2 \aH^n_k(B^m_{2R_2}(\bar{x}_k)),
        \]
        from which our claim \eqref{E24} follows. \par 
        Finally, using the mean value property for sub-harmonic function on the graph of each $u_k$ (i.e. \eqref{Assumption}) we get,
        thanks to \eqref{E22} and  \eqref{E24}, that
        \begin{align*}
        \sup_{B_{2}(x_k)} |\nabla u_k|^2 & \leq \sup_{B^m_{R_2/2}(x_k)} |\nabla u_k|^2 \leq \fint_{B^m_{R_2}(x_k)} |\nabla u_k|^2 \, d\aH^n_k
        \\
        & \leq (1+4c^2)\fint_{B^m_{R_2}(x_k)} |\nabla_m u_k|^2 \, d\aH^n_k \leq 4 \epsilon^2 (1+4c^2)
        \frac{\aH^n_k(B^m_{2R_2}(\bar{x}_k))}{\aH^n_k(B^m_{R_2}(\bar{x}_k))} \leq C(u,n) \epsilon^2.
        \end{align*}
        In particular, by Lemma \ref{L|uniform lip}, we get
        \[
        \sup_{B_{1}(x)} |\nabla u|^2 \leq C\epsilon^2.
        \]
        Since $\epsilon>0$ was arbitrary we conclude.
    \end{proof}
\end{customthm}

\printbibliography
\end{document}